\documentclass[11pt, reqno]{amsart}

\usepackage{amssymb}
\usepackage{hyperref} 
\usepackage{mathrsfs,bbold}
\usepackage{amsmath}
\usepackage{amsthm}
\usepackage{enumerate}
\usepackage{dsfont}
\usepackage{color}
\usepackage[a4paper]{geometry}
\usepackage[utf8]{inputenc}
\usepackage{stmaryrd}
\usepackage{titlesec}
\usepackage{mathabx}
\usepackage{appendix}
\usepackage{todonotes}

\usepackage[all]{xy}
\let\oldtocsection=\tocsection
\let\oldtocsubsection=\tocsubsection
\let\oldtocsubsubsection=\tocsubsubsection

\renewcommand{\tocsection}[2]{\hspace{0em}\oldtocsection{#1}{#2}}
\renewcommand{\tocsubsection}[2]{\hspace{1em}\oldtocsubsection{#1}{#2}}
\renewcommand{\tocsubsubsection}[2]{\hspace{2em}\oldtocsubsubsection{#1}{#2}}

\numberwithin{theorem}{section}
\numberwithin{equation}{section}

  \newcommand{\ZZ}{\mathbb{Z}}    
\newcommand{\RR}{\mathbb{R}}  \newcommand{\CC}{\mathbb{C}}      
   \newcommand{\VV}{\mathbb{V}}

 \newcommand{\mcC}{\mathcal{C}}  \newcommand{\mcF}{\mathcal{F}}

   \newcommand{\fF}{{\frak{F}}}

\newcommand{\bfV}{{\bf V}}    \newcommand{\xx}{\mathbf{x}}

\newcommand{\fA}{\mathfrak{A}}
\newcommand{\fB}{\mathfrak{B}}

\newcommand{\fM}{\mathfrak{M}}

\newcommand{\fV}{\mathfrak{V}}

\newcommand{\fX}{\mathfrak{X}}

\newcommand{\EE}{\mathbb{E}} \newcommand{\PP}{\mathbb{P}}

\newcommand{\eps}{\epsilon}
\renewcommand{\leq}{\leqslant}
\renewcommand{\geq}{\geqslant}

\newcommand{\ssk}{\smallskip}

\newcommand{\cA}{\mathcal{A}}
\newcommand{\cB}{\mathcal{B}}
\newcommand{\cC}{\mathcal{C}}

\newcommand{\cF}{\mathcal{F}}
\newcommand{\cG}{\mathcal{G}}

\newcommand{\cL}{\mathcal{L}}

\newcommand{\bV}{\mathbf{V}}

\newcommand{\dd}{\mathrm{d}}

\newenvironment{Dem}[1][\unskip]{%
    \begin{list}{\hspace{0.5cm}{\sf \textbf{Proof #1 --}}}{%
        \setlength{\topsep}{0pt}%
        \setlength{\leftmargin}{0pt}%
        \setlength{\rightmargin}{0pt}%
        \setlength{\listparindent}{0pt}%
        \setlength{\itemindent}{0pt}%
        \setlength{\parsep}{0pt}%
        \addtolength{\leftmargin}{20pt}%
        \addtolength{\rightmargin}{0pt}%
    } \item }{\hfill{\begin{flushright} $\rhd$ \end{flushright}}\end{list}\smallskip}

\newenvironment{DemThm}{%
    \begin{list}{\hspace{0.5cm}{\textbf{\textsf{Proof} --}}}{%
        \setlength{\topsep}{0pt}%
        \setlength{\leftmargin}{0pt}%
        \setlength{\rightmargin}{0pt}%
        \setlength{\listparindent}{0pt}%
        \setlength{\itemindent}{0pt}%
        \setlength{\parsep}{0pt}%
        \addtolength{\leftmargin}{20pt}%
        \addtolength{\rightmargin}{0pt}%
    } \item }{\hfill{\space $\rhd$}\end{list}\smallskip}

\titleformat{\section}[block]
{\filcenter\normalfont\sffamily\bfseries\Large}{{\hspace{-0.7cm}}\thesection \vspace{0.3cm}}{0.75em}{}

\titleformat{\subsection}[block]
{\normalfont\sffamily\bfseries\large}{\hspace{-1cm}\thesubsection \vspace{0.3cm}}{.75em}{}
\titlespacing{\subsection}{-0pc}{1.5ex plus .1ex minus .2ex}{0pc}

\titleformat{\subsubsection}[block]
{\normalfont\sffamily\bfseries}{\hspace{-1.2cm}\thesubsubsection  \vspace{0.3cm}}{.75em}{}
\titlespacing{\subsection}{0pc}{1.5ex plus .1ex minus .2ex}{0pc}

\def\XXint#1#2#3{{\setbox0=\hbox{$#1{#2#3}{\int}$}
     \vcenter{\hbox{$#2#3$}}\kern-.5\wd0}}

\numberwithin{subsection}{section}
\numberwithin{subsubsection}{subsection}

\newtheoremstyle{mystyle}
{3pt}
{3pt}
{\sffamily}
{0em}
{\bfseries\sffamily}
{.}
{.5em}
{}
\theoremstyle{mystyle}

\newtheorem{thm}{Theorem}
\newtheorem*{thm*}{Theorem}

\newtheorem{lem}[thm]{\hspace{-0.15cm}  {Lemma} }
\newtheorem{prop}[thm]{\hspace{-0.15cm} {Proposition}}
\newtheorem{defn}[thm]{ \hspace{-0.3cm} {Definition}}
\newtheorem*{defn*} {Definition}
\newtheorem*{prop*} {Proposition}
\newtheorem*{lem*} {Lemma}
\newtheorem*{cor*} {Corollary}

\numberwithin{equation}{section} 

\author{R\'emi Catellier}
\email{remi.catellier@univ-rennes1.fr} 
\address{Universit\'e de Rennes 1 - IRMAR - Centre Henri Lebesgue}
\author{Isma\"el Bailleul}
\email{ismael.bailleul@univ-rennes1.fr}
\address{Universit\'e de Rennes 1 - IRMAR}
\date{\today}
\title{From rough flows to homogenization for fast-slow dynamics}
\keywords{Rough flows, random ordinary differential equations, stochastic homogenization, invariance principle}

\begin{document}

\vspace*{3ex minus 1ex}
\begin{center}
\huge\sffamily{Rough flows and homogenization in stochastic turbulence}
\end{center}
\vskip 5ex minus 1ex

\begin{center}
{\sf I. BAILLEUL\footnote{I.B. thanks the U.B.O. for their hospitality.} and R. CATELLIER\footnote{R. Catellier is supported by the Labex Lebesgue  \\  
{\sf AMSClassification:} 60F05, 37H10}}
\end{center}

\vspace{1cm}

\begin{center}
\begin{minipage}{0.8\textwidth}
\renewcommand\baselinestretch{0.7} \rmfamily {\scriptsize {\bf \sc \noindent Abstract.} We provide in this work a tool-kit for the study of homogenisation of random ordinary differential equations, under the form of a friendly-user black box based on the tehcnology of rough flows.  We illustrate the use of this setting on the example of stochastic turbulence. 
}
\end{minipage}
\end{center}

\bigskip

\setcounter{tocdepth}{2}
\tableofcontents

\section[\hspace{0.6cm} {\sf Introduction}]{Introduction}
\label{SectionIntroduction}

The history of averaging and homogenization problems for dynamical systems is fairly long and has its roots in classical perturbative problems in mechanics, in the 19th centery. It has evolved in an impressive body of methods and tools used to analyse a whole range of multiscale systems, such as (possibly random) transport equations with multiple time-scales \cite{PavliotisStuart, Kabanov}, or heat propagation in random media \cite{CioranescuDonato, Tartar}. The latest developments of Otto, Gloria \& co \cite{Otto} and Armstrong \& co \cite{ArmstrongKuusiMourrat, ArmstrongGloriaKuusi} on homogenization for the solutions of Hamilton-Jacobi equations use and develop deep results in partial differential equations. The present work deals with the transport side of the story, in the line of the classical works of Kesten and Papanicolaou on homogenization for random stochastic differential equations \cite{KestenPapanicolaou, KP2, PSV}, and put them in the flow of ideas and tools that have emerged in the early 2000's with rough paths theory. Kelly and Melbourne \cite{KellyMelbourne1, KellyMelbourne2} have for instance shown recently how one can use rough paths methods to investigate a fast-slow system of the form 
$$
\dot x_\epsilon = a(x_\epsilon,y_\epsilon) + \frac{1}{\epsilon}\,b(x_\epsilon,y_\epsilon),
$$
where the dynamics of the fast component $y_\epsilon$ is autonomous and Anosov or axiom A, or even non-uniformly hyperbolic. We would like to put this result and other homogenization results in the newly introduced setting of \textit{rough flows} \cite{BailleulRiedel}, that encompasses a large part of the theory of rough differential equations, and unifies it with the theory of stochastic flows. We provide for that purpose an easily usable black box for the study of homogenisation of random ordinary differential equations, under the form of a result   \vspace{0.15cm}

\begin{center}
\textit{Convergence of finite dimensional marginals} $\oplus$ \textit{Moment/tightness bounds}   \vspace{0.1cm}

{\small (for the driving vector fields)}   \vspace{0.15cm}

$\Longrightarrow$ \textit{Homogenisation}
\end{center}

\noindent for which no knowledge of the mechanics of rough flows is required. See Theorem \ref{ThmUseIt} in section \ref{SubsectionUseIt}. As an illustration of use of this method in homogenization problems, we show in the present work how one can get back and extend in a clean and efficient way Kesten and Papanicolaou' seminal result \cite{KestenPapanicolaou} on stochastic turbulence.

\ssk

The theory of rough flows is based on the following paradigm. The kind of dynamics we are about to consider are all generated by some kind of time-dependent vector fields, or drivers, that generate flows by a \textit{deterministic continuous} mechanism. Any ordinary differential equation is naturally recast in this setting. The benefits of this picture for the study of averaging and homogenization problems are obvious. If the drivers are random and depend on some parameters, it suffices that they converge in law in the space of drivers for their associated dynamics to converge in law, from the continuity of the driver-to-flow map. Support theorems and large deviation results are also automatically transported from the driver world to the flow world. The rough flow setting somehow provides an optimized and friendly environment where to apply ideas similar to those of rough paths theory, with the same benefits. As a matter of fact, one can also study some homogenisation problems for random ordinary differential equations from the latter point of view, such as done by Kelly and Melbourne in their works \cite{KellyMelbourne1, KellyMelbourne2} on fast-slow systems, to the expense of working with tensor products of Banach spaces and the involved subtleties. No such high level technology is required in the elementary setting of rough drivers and rough flows, which may then be easier to use \cite{BailleulRiedel}. More importantly, it has a dual version on function spaces that can be used to study some hyperbolic partial differential equations and seem beyond the scope of Lyons' formulation of rough paths theory \cite{BailleulGubinelli}.

\medskip

Section \ref{SubsectionBlackBox} provides a very light presentation of rough drivers and their associated flows; convergence problems for flows amount in this setting to convergence problems for their drivers -- a philosophy shared by the martingale problem formulation of stochastic dynamics, with the noticeable difference that we are here in a deterministic setting. Section \ref{SubsectionUseIt} contains the above generic homogenisation result; it is proved in Appendix \ref{SectionAppendixCompactness}. An elementary deterministic example is given in section \ref{SubsectionPureArea} as an illustration of the mechanics at play in the rough driver/flow setting. The case of homogenisation for stochastic turbulence is treated in section \ref{SectionKestenPapanicolaou}.

\bigskip

\noindent {\bf Notations.} We shall use the sign $\lesssim$ for an inequality that holds up to a multiplicative positive constant whose precise value is unimportant. The sign $\lesssim_c$ will be used to indicate that this constant depends on a parameter $c$. Given a finite positive time horizon $T$, we shall write $D_T$ for $\big\{(s,t)\in[0,T]^2\,;s\leq t\big\}$. We shall use the $|\cdot|$ sign to denote any Euclidean norm on a finite dimensional space; its precise choice will be unimportant.

\ssk

\begin{itemize}
   \item Given a non-integer positive regularity index $a$, we shall denote by $\mcC^a$ the usual space of $a$-H\"older functions. Given $0<a_1<1$, a $2$-index map $(Z_{ts})_{0\leq s\leq t\leq T}$ with values in some space $\mcC^{a_2}_b(\RR^d)$ will be said to be $a_1$-H\"older if 
$$
\big\|Z\big\|_{\mcC^{a_1}_{ts} \mcC^{a_2}_b} := \underset{0\leq s<t\leq T}{\sup}\;\frac{\big\|Z_{ts}\big\|_{\mcC^{a_2}_b}}{|t-s|^{a_1}} < \infty;
$$
we write $Z\in \mcC^{a_1}_{ts} \mcC^{a_2}_b$.   \vspace{0.1cm}
   
   \item An \textit{additive function of time} $(V_{ts})_{0\leq s\leq t\leq T}$ is a vector space valued function $V$ of time such that $V_{ts} = V_{tu} + V_{us}$, for all $0\leq s\leq u\leq t\leq T$.   \vspace{0.1cm}
   
   \item Whenever convenient, we shall freely identify vector fields with first order differential operators, so that given two vector fields $V_1, V_2$, the notation $V_1V_2$ will stand for the second order differential operator whose action on smooth functions $f$ is 
   $$
   V_1V_2f = (Df)\big(DV_2(V_1)\big) + (D^2f)(V_1,V_2).   
   $$  \vspace{0.1cm}
   
   \item Given $f\in L^\infty(\RR^d,\RR^d)$  and $\sigma$ in the unit ball of $\RR^d$, we define inductively a sequence $\Delta_\sigma^m$ of operators on $L^\infty(\RR^d,\RR^d)$ setting
$$
\big(\Delta_\sigma f\big)(\cdot) = f(\cdot + \sigma) - f(x) \quad \mathrm{and} \quad \Delta^{m+1}_\sigma f = \Delta_\sigma (\Delta^m_\sigma f).
$$   \vspace{0.1cm}

   \item Implicit summation of repeated indices is used throughout, so $a^ib_i$ means $\sum_i a^i b_i$.   \vspace{0.1cm}
   
   \item We denote by $L^k(\Omega)$ the corresponding integrability spaces over some probability space $(\Omega,\mcF,\PP)$. 
\end{itemize}

\bigskip

\section[\hspace{0.6cm} {\sf Tools for flows of random ODEs}]{Tools for flows of random ODEs}
\label{SectionToolsFlows}

The machinery of rough drivers and rough flows introduced in \cite{BailleulRiedel} provides a very convenient setting for the study of convergence of flows and weak convergence of random flows. Rather than giving the reader an account of the theory of rough flows, we single out here part of it under the form of a friendly user black box that requires no knowledge of the mechanics of rough flows. We refer the interested reader to the work \cite{BailleulRiedel}, and to Appendix \ref{SectionAppendixCompactness}, for some more technical details.

\bigskip

\subsection[\hspace{-0.6cm} {\sf The black box}]{The black box.}
\label{SubsectionBlackBox}

The starting point of this business is the elementary observation that if we are given some smooth globally Lipschitz vector fields $v_1,\dots,,v_\ell$ on $\RR^d$, and some real-valued controls $h^1,\dots,h^\ell$ on some time interval $[0,T]$, then the solution flow $(\varphi_{ts})_{0\leq s\leq t\leq T}$ of the controlled ordinary differential equation
$$
\dot z_t = \dot h^i_t \, v_i(z_t)
$$
enjoys the Taylor expansion property
\begin{equation}
\label{EqTaylorFormula}
f\circ \varphi_{ts} = f + \big(h^i_t-h^i_s\big)V_if  + \left(\int_s^t\int_s^{u_1} dh^j_{u_2}\,dh^k_{u_1}\right)\,V_jV_kf + O\big(|t-s|^{>2}\big)
\end{equation}
for all smooth functions $f$. The notion of rough driver captures the essence of the different terms that appear in this local description of the dynamics.

\medskip

\begin{defn}   {\it
Let $2\leq p<2+r<3$ be given. A \textbf{\textsf{\emph{rough driver}}}, with regularity indices $p$ and $(2+r)$, is a family $\big({\bfV}_{ts}\big)_{0\leq s\leq t\leq T}$, with ${\bfV}_{ts} := \big(V_{ts},\VV_{ts}\big)$, for some vector fields $V_{ts}$, and $\VV_{ts}$ some second order differential operator, such that \vspace{0.1cm}
\begin{itemize}
   \item[{\bf (i)}] the vector field $V_{ts}$ is an additive function of time, with $V\in \mcC^{1/p}_{ts} \mcC^{2+r}_b$,   \vspace{0.15cm}

   \item[{\bf (ii)}] the second order differential operators 
\begin{equation*}
W_{ts} := \VV_{ts} - \frac{1}{2}V_{ts}V_{ts}, 
\end{equation*}
   are actually vector fields, and $W\in \mcC^{2/p}_{ts} \mcC^{1+r}_b$,   \vspace{0.15cm}

   \item[{\bf (iii)}] we have 
   $$
   \VV_{ts} = \VV_{tu} + V_{us}V_{tu} + \VV_{us}, 
   $$
   for any $0\leq s\leq u\leq t\leq T$.   
\end{itemize}
We define the norm of $\bf V$ to be 
\begin{equation*}
\|{\bfV}\| := \max\left(\big\|V\big\|_{\mcC^{1/p}_{ts}\mcC^{2+r}_b}\,,\, \big\|W\big\|_{\mcC^{2/p}_{ts}\mcC^{1+r}_b}\right).
\end{equation*}   }
\end{defn}

\medskip

We simply talk of a rough driver when its regularity indices are clear from the context. We typically use rough drivers to give a local description of the dynamics of a flow $\varphi$, under the form of a Taylor expansion formula
$$
f\circ\varphi_{ts} \simeq f + V_{ts}f + \VV_{ts}f.
$$
In the Taylor formula \eqref{EqTaylorFormula}, the term $\big(h^i_t-h^i_s\big)V_i$ plays the role of $V_{ts}$, while the term $\left(\int_s^t\int_s^r dh^j_u\,dh^k_r\right)\,V_jV_k$ has the role of $\VV_{ts}$; check that properties {\sf (i)-(iii)} hold indeed for these two terms. More generally, given any sufficiently regular time-dependent vector field $v_t$ on $\RR^d$, on can check that setting
\begin{equation}
\label{EqLiftSmoothVectorField}
V_{ts} := \int_s^t v_u \, du,\quad \VV_{ts} := \int_s^t\int_s^{u_1} v_{u_2} v_{u_1}\, du_2 du_1
\end{equation}
defines a rough driver for which
$$
W_{ts} = \int_s^t\int_s^r \big[v_{u_2}, v_{u_1}\big]\, du_2 du_1,
$$
with Lie brackets of vector fields used here. Formula \eqref{EqLiftSmoothVectorField} defines the \textsf{\textbf{canonical lift}} of a possibly time-dependent vector field $v$. As we shall use it later, remark here that if ${\bfV} = \big(V, W+\frac{1}{2}\,V^2\big)$ stands for a rough driver with regularity indices $p$ and $(2+r)$, and $X$ stands for a $\frac{2}{p}$-H\"older function with values in the space of $\mcC^{1+r}_b$ vector fields on $\RR^d$, then the formula 
$$
\big(V_{ts}, W_{ts}+\frac{1}{2}\,V_{ts}^2 + X_t-X_s\big)
$$
still defines a $(p,2+r)$-rough driver.

\medskip

\begin{defn}   {\it
Let $V_0$ be a bounded Lipschitz vector field on $\RR^d$; let also $\bf V$ be a rough driver with regularity indices $p$ and $(2+r)$. A \textsf{\emph{\textbf{flow}}} $\big(\varphi_{ts}\big)_{0\leq s\leq t\leq T}$ is said to \textbf{\textsf{\emph{solve the rough differential equation}}}
\begin{equation*}
d\varphi = V_0(\varphi)dt + {\bfV}(\varphi\,; dt)
\end{equation*}
if there exists a possibly $\big(V_0,\bfV\big)$-dependent positive constant $\delta$ such that the inequality
\begin{equation*}
\Big\|f\circ\varphi_{ts} - \Big\{f + (t-s)\big(V_0f\big) + V_{ts}f + \VV_{ts}f \Big\}\Big\|_\infty \lesssim  \|f\|_{\mcC^{2+r}}|t-s|^\frac{3}{p}
\end{equation*}
holds for all $f\in\mcC^{2+r}_b$, and all $0\leq s\leq t\leq T$ with $t-s\leq\delta$. Such flows are called \textbf{\textsf{\emph{rough flows}}}.   }
\end{defn} 

\medskip

If $\bfV$ is the canonical lift of a $\mcC^{2+r}_b$ time-dependent vector field $v$, its associated rough flow coincides with the classical flow generated by $v$. A robust well-posedness result is provided by the next result, proved in \cite{BailleulRiedel}.

\medskip

\begin{thm} 
\label{ThmConstructingRoughFlows}   {\it 
Assume $\frac{p}{3}<r\leq 1$. Then the differential equation on flows
$$
d\varphi = V_0(\varphi)dt + {\bfV}(\varphi \,; dt)
$$ 
has a unique solution flow; it takes values in the space of homeomorphisms of $\RR^d$, and depends continuously on $V_0$ and $\bfV$ in the topology of uniform convergence. Moreover, if $r<1$, then the maps $\varphi_{ts}$ and their inverse  have uniformly bounded $\mcC^r$-norms; if $r=1$, they have uniformly bounded Lipschitz norms.   }
\end{thm} 

\medskip

If $B$ is an $\ell$-dimensional Brownian motion and $v_1, \dots, v_\ell$ are $\mcC^3_b$ vector fields on $\RR^d$, one can prove that setting
$$
V_{ts} = B^i_{ts} v_i, \qquad \textrm{and}  \qquad \VV_{ts} = \left(\int_s^t\int_s^{u_1} \circ dB^j_{u_2}\circ dB^k_{u_1}\right) v_jv_k
$$
defines almost surely a rough driver with regularity indices $p$ and $(2+r)$, for any $p<2+r<3$ with $\frac{p}{3}<r$, and that the solution flow of the equation
$$
d\varphi = {\bfV}(\varphi\,; dt)
$$
coincides almost surely with the flow generated by the Stratonovich stochastic differential equation
$$
dx_t = v_i(x_t)\,{\circ dB^i_t}.
$$
See e.g. Lyons' seminal paper \cite{Lyons98}.

\bigskip

\subsection[\hspace{-0.6cm} {\sf How to use it}]{How to use it.}
\label{SubsectionUseIt}

Let then assume we are given a random ordinary differential equation
\begin{equation}
\label{EqRandomODE}
\dot x^\epsilon_t = v^\epsilon_t\big(x^\epsilon_t\big)
\end{equation}
in $\RR^d$, driven by a random time-dependent globally Lipschitz vector field $v^\epsilon_t$, depending on a parameter $\epsilon$, an element of $(0,1]$ say. One can think for instance of the \textit{slow dynamics in a fast-slow system} \cite{KellyMelbourne2}
\begin{equation*}
\begin{split}
& \dot x^\epsilon_t = f\big(x^\epsilon_t,y^\epsilon_t\big),   \\
& \dot y^\epsilon_t = \frac{1}{\epsilon}\, g\big(y^\epsilon_t\big),
\end{split}
\end{equation*}
driven by deterministic vector fields $f,g$, but where $y^\epsilon_0$ is random for instance, so we have \eqref{EqRandomODE} with
$$
v^\epsilon_t(\cdot) =  f\big(\cdot, y^\epsilon_t\big).
$$ 
We shall also reformulate in section \ref{SectionKestenPapanicolaou} the \textit{stochastic turbulence} dynamics in those terms. Fix a finite time horizon $T$ and define, for $0\leq s\leq t\leq T$, the canonical lift ${\bfV}^\epsilon$ of $v^\epsilon$ into a rough driver
$$
V^\epsilon_{ts} := \int_s^t v^\epsilon_u \, du,\quad \VV^\epsilon_{ts} := \int_s^t\int_s^{u_1} v^\epsilon_{u_2} v^\epsilon_{u_1} \, du_2 du_1,
$$
and
$$
W^\epsilon_{ts} := \int_s^t\int_s^r \big[v^\epsilon_{u_2}, v^\epsilon_{u_1}\big]\, du_2 du_1.
$$
Denote by $\varphi^\epsilon$ the random flow generated by equation \eqref{EqRandomODE}, so $\varphi^\epsilon_{ts}(x)$ is, for any $0\leq s\leq t\leq T$, the value at time $t$ of the solution to equation \eqref{EqRandomODE} started from $x$ at time $s$. This flow is also the solution flow of the equation 
$$
d\varphi^\epsilon = {\bfV}^\epsilon\big(\varphi^\epsilon\,;dt\big).
$$ 
Given $0<r<1$, denote by $\mcC^r_{(0)}$ the space of $r$-H\"older continuous functions from $\RR^d$ to itself that are at finite $\mcC^r$-distance from the identity, and write $\textsf{Diff}^r_{(0)}$ for the space of $\mcC^r$-homeomorphisms with $\mcC^r$-inverse, for which both the homeomorphism and its inverse are at finite $\mcC^r$-distance from the identity. The following convergence result, proved in Appendix \ref{SectionAppendixCompactness}, is an elementary   \vspace{0.15cm}

\begin{center}
\textit{Convergence of finite dimensional marginals} $\oplus$ \textit{Moment/tightness bounds}   \vspace{0.1cm}

{\small (for the driving vector fields)}   \vspace{0.15cm}

$\Longrightarrow$ \textit{Homogenisation}
\end{center}

\noindent result. The exponent $a$ in the statement is to be thought of as a big positive constant. 

\medskip

\begin{thm}
\label{ThmUseIt}   {\it 
Let some positive finite exponents $(p,r\,;a)$, with $r<1$, be given such that 
$$
0 < \frac{1}{\frac{1}{p}-\frac{1}{2a}} - 2 < r - \frac{d}{a}.
$$
\begin{itemize}
   \item[\textcolor{gray}{$\bullet$}] Assume that for each $0\leq s\leq t\leq T$, and each $y\in\RR^d$, the random variables $V^\epsilon_{ts}(y)$ and $W^\epsilon_{ts}(y)$ converge weakly as $\epsilon$ goes to $0$.   \vspace{0.15cm}
   
   \item[\textcolor{gray}{$\bullet$}] Assume further that there is an integer $k_1\geq 3$ for which the positive quantity
   {\small \begin{equation*}
   \begin{split}
   \int_{\RR^d} \left\{\left\|\frac{V^\epsilon_{ts}(y)}{|t-s|^\frac{1}{p}}\right\|^a_{L^{2a}(\Omega)} + \left\|\frac{W^\epsilon_{ts}(y)}{|t-s|^\frac{2}{p}}\right\|^a_{L^a(\Omega)} \right\}\,dy &+ \iint_{\RR^d\times B(0,1)}  \left\|\frac{\Delta_\sigma^{k_1} V^\epsilon_{ts}(y)}{|t-s|^\frac{1}{p}}\right\|^a_{L^{2a}(\Omega)}\,dy\,\frac{d\sigma}{|\sigma|^{(2+r)a+d}}   \\
   &+ \iint_{\RR^d\times B(0,1)} \left\|\frac{\Delta_\sigma^{k_1-1} W^\epsilon_{ts}(y)}{|t-s|^\frac{2}{p}}\right\|^a_{L^a(\Omega)}\,dy\,\frac{d\sigma}{|\sigma|^{(1+r)a+d}}
   \end{split}
   \end{equation*}   }
   is bounded above by a finite constant independent of $\epsilon$.   \vspace{0.1cm}
\end{itemize}

Then for every pair of regularity indices $(p',2+r')$, with $p'<2+r'<3$, and
$$
r'<r-\frac{d}{a}, \qquad \textrm{and}  \qquad \frac{1}{3} < \frac{1}{p'} < \frac{1}{p} - \frac{1}{2a},
$$
there exists a random rough driver $\bfV$, with regularity indices $p'$ and $(2+r')$, whose associated random flow 
$$
d\varphi = {\bf V}(\varphi\,;dt)
$$
is the weak limit in $C\Big([0,T],\sf{Diff}^{r'}_{(0)}\Big)$ of the random flows $\varphi^\epsilon$ generated by the dynamics \eqref{EqRandomODE}.   }
\end{thm}

\medskip

The above $\epsilon$-uniform moment bound is actually a sufficient condition for \textit{tightness} in the space of drivers with regularity indices $p'$ and $(2+r')$. Note that the convergence and moment assumptions are about the vector fields $V^\epsilon$ and $W^\epsilon$ that generate the dynamics, while the conclusion is on the dynamics itself. The possibility to transfer a weak convergence result on the rough drivers to the dynamics comes from the continuity of the solution map, given as a conclusions in Theorem \ref{ThmConstructingRoughFlows}. Note also that we work here with vector fields $V,W$ that are in particular bounded, as required by the definition of a $\mcC^\alpha$ function, for a non-integer regularity exponent $\alpha$. In applications, one may have first to localize the dynamics in a big ball of radius $R$, use Theorem \ref{ThmUseIt}, and remove the localization in a second step. This is what we shall be doing in our study of stochastic turbulence in section \ref{SectionKestenPapanicolaou}. Let us note here that the results proved by Kelly and Melbourne \cite{KellyMelbourne1,KellyMelbourne2} in their study of fast-slow systems with a chaotic fast component can actually be rephrased exactly in the terms of Theorem \ref{ThmUseIt}, so one can get back their conclusions from the point of view developed here.

\bigskip

\subsection[\hspace{-0.6cm} {\sf A toy example}]{A toy example.}
\label{SubsectionPureArea}

Before applying Theorem \ref{ThmUseIt} in the setting of stochastic turbulence, we illustrate in this section on an elementary and interesting toy example the fundamental continuity property of the solution map to an equation
$$
d\varphi = {\bfV}(\varphi\,;dt),
$$
in a deterministic setting. In this example, we construct a family ${\bfV}^\epsilon$ of rough drivers, obtained as the canonical lift of a smooth $\epsilon$-dependent vector field on the plane, such that its first level $V^\epsilon$ converges to $0$ in a strong sense while the flow $\varphi^\epsilon$ associated with ${\bfV}^\epsilon$ does not converge to the identity. This shows the crucial influence of the second level object $\VV^\epsilon$ on the dynamics generated by ${\bfV}^\epsilon$. We work in $\RR^2\simeq \CC$ ; set
$$
v_t(x) := if(x)\,e^{if(x)t}, 
$$
for some $\mcC^3_b$ non-zero phase $f$, so that its canonical lift ${\bfV} = \big(V,\frac{1}{2}\,V^2 + W\big)$ as a rough driver has first level
$$
V_{ts}(x) = e^{if(x)t} - e^{if(x)s}.
$$ 
Given $2\leq p<3$, we define a space/time rescaled rough driver ${\bfV}^\epsilon$, with regularity indices $p$ and $1$, setting
$$
V^\epsilon_{ts}(x) := \epsilon\,V_{t\epsilon^{-2}\,s\epsilon^{-2}}\big(\epsilon^2\,x\big), \quad W^\epsilon_{ts}(x) := \epsilon^4\,\Big(W_{t\epsilon^{-2}\,s\epsilon^{-2}}\Big)\big(\epsilon^2\,x\big);
$$ 
this is the canonical lift of the $\epsilon$-dependent vector field
$$
v^\epsilon_t(x) := \frac{1}{\epsilon}\,v(\epsilon^2x) = \frac{i\,f(\epsilon^2 x)}{\epsilon}\,e^{if(\epsilon^2 x)t}.
$$

\medskip

\begin{thm}\label{theorem:pure_area}   {\it 
The rough driver ${\bfV}^\epsilon$ converges as a $(p,1)$-rough drivers to the pure second level rough driver
$$
{\bfV}_{t,s}(x) := \left(0 , -\frac14 \big(t^2 - s^2\big)\;f(0)\,(\nabla f)(0)\right).
$$   }
\end{thm}

\medskip

As a corollary, the solution flow $\varphi^\epsilon$ to the equation 
$$ 
\dot x^\epsilon_t = v^\epsilon_t(x^\epsilon_t)
$$
converges to the elementary flow generated by the ordinary differential equation
$$
\dot x_t = -\frac{1}{2}\,f(0)\,(\nabla f)(0)
$$
with constant vector field.

\medskip
 
\begin{DemThm}
We shall prove the claim as a direct consequence of the following elementary estimate
\begin{equation}
\label{EqElementaryEstimate}
\Big\| D^\ell\Big(e^{if(\cdot)t} - e^{if(\cdot)s}\Big)\Big\|_\infty \leq_{f,\gamma} T_1^\ell\,|t-s|^\gamma,
\end{equation}
that holds for all times $0\leq s\leq t\leq T_1<\infty$, every exponent $0<\gamma\leq 1$, and any derivative index $0\leq\ell\leq 3$, as shown by interpolating two trivial bounds.

\ssk

Working with $T_1=T\epsilon^{-2}$, and since $D^\ell V^\eps_{t,s}(x) = \eps^{2\ell+1} \Big(D^\ell V^\eps_{t\eps^{-2},s\eps^{-2}}\Big)\big(x\eps^2\big)$, it already follows from \eqref{EqElementaryEstimate} that 
$$
\Big\| D^\ell V^\epsilon_{ts} \Big\|_\infty \lesssim T^\ell\epsilon^{1-2\gamma}\,|t-s|^\gamma,
$$
so, indeed, we have 
$$
\sup_{0\le s \le t \le T} \frac{\|V^\eps_{t,s}\|_{C^3}}{|t-s|^\gamma} \longrightarrow_{\eps \to 0} 0
$$
if one chooses $0<\gamma<\frac{1}{2}$. 

\ssk

To deal with $W^\epsilon$, note first that an integration by parts gives for $W$ the decomposition
{\small \begin{align*}
W_{t,s}(x) :=& \frac12 \int_s^t \Big(D_xv_{r}\,V_{rs}(x) - D_xV_{r,s}\,v_r(x)\Big)\,dr  = \frac12 DV_{ts}(x) V_{t,s}(x) - \int_s^t DV_{rs}(x) v_r(x) \,dr   \\
                   =& \frac12  DV_{ts}(x) V_{t,s}(x)  -  \int_s^t if(x)\Big(r e^{if(x) r} - s e^{if(x)s}\Big) \underbrace{\frac12\Big(\nabla f(x) \overline{ie^{if(x) r}} + \overline{\nabla f(x)} i e^{if(x) r}\Big)}_{\langle \nabla f(x), v_r(x) \rangle}\,dr;
\end{align*}   }
so one can write
$$
W_{t,s}(x) = -\frac14 (t^2-s^2) f(x) \nabla f(x) + R_{ts}(x)
$$
with
\begin{align*}
 R_{t,s}(x) = &\frac12 DV_{ts}(x) V_{ts}(x) + \frac s 2 e^{i f(x) s} f(x) \int_s^t \Big( \nabla f(x) - \overline{\nabla f(x)} e^{2 i f(x) r}\Big) \,dr \\
 					& + f(x) \bar{\nabla f(x)} \frac 12 \int_s^t r e^{2i f(x) r } \,dr.
\end{align*}
Hence 
\begin{align*}
W^\eps_{t,s}(x) 
&= \eps^4 W_{t\eps^{-2},s\eps^{-2}}(x\eps^2)\\
&=  -\frac14 (t^2-s^2) b(\eps^2 x) \nabla b(\eps^2x) + R^\eps_{ts}(x),
\end{align*}
where
$$
R^\eps_{ts}(x) := \eps^4 R_{t\eps^{-2}\,s\eps^{-2}}(x\eps^2).
$$

The scaling in $\epsilon$ between space and time gives to convergence
$$
\sup_{0 \le s \le t \le T} \frac{\| R^\eps_{t,s}\|_{C^2}}{|t-s|^{2\gamma}}  \to 0,
$$
which is enough to conclude that $W^\eps$ converges in the same space to 
$$
(t,s,x) \to -\frac14 (t^2 - s^2) f(0) \nabla f (0).
$$  \vspace{-0.3cm}
\end{DemThm}

\bigskip

\section[\hspace{0.6cm} \textsf{A case study: Stochastic turbulence}]{A case study: Stochastic turbulence}
\label{SectionKestenPapanicolaou}

We show in this section how one can use the black box provided by Theorem \ref{ThmUseIt} to reprove and improve in a simple way Kesten and Papanicolaou' seminal result on stochastic turbulence \cite{KestenPapanicolaou}. The object of interest here is the dynamics of a particle subject to a random velocity field that is a small perturbation of a constant deterministic velocity. Precisely, consider the random ordinary differential equation
$$
\dot x_t = {\sf v} + \epsilon F(x_t),
$$
with initial condition $x_0$ fixed, where {\sf v} is a deterministic non-zero mean velocity and $F$ is a sufficiently regular centered, stationary, random field; precise assumptions are given below. To investigate the fluctuations of $x_\bullet$ around its typical value, one looks at the dynamics of the recentered and time-rescaled process
$$
x^\epsilon_t := x_{\epsilon^{-2} t} - \epsilon^{-2} t\,{\sf v},
$$ 
and prove that the continuous random processes $\big(x^\epsilon_t\big)_{0\leq t\leq 1}$ converge in law, as $\epsilon$ decreases to $0$, to a Brownian motion with some constant drift $b$ and some covariance $\sigma^*\sigma$, both given explicitly in terms of the statistics of $F$. We actually use Theorem \ref{ThmUseIt} to prove a similar result for flows directly.

\medskip

As the process $x^\epsilon$ solves the random ordinary differential equation
\begin{equation}
\label{EqStochasticTurbulenceEq}
\dot x^\epsilon_t = \epsilon^{-1}F\big(x^\epsilon_t + \epsilon^{-2}t{\sf v}\big),
\end{equation}
the flow genetared by the latter dynamics is also associated with the "rough driver"
\begin{equation}
\label{EqDefnVep}
{\bfV}^\epsilon = \Big(V^\epsilon,W^\epsilon+\frac{1}{2}\,(V^\epsilon)^2\Big),
\end{equation}
where
$$
V^\eps_{t,s}(x) := \frac1\eps\int_s^t F\Big(x+ \frac{u}{\epsilon^2}\,{\sf v}\Big)\,du  =: \int_s^t v^\epsilon_u(x)\,du
$$
and 
$$
W^\eps_{t,s}(x) := \frac{1}{2} \, \int_s^t\int_s^{u_1} \big[V^\epsilon_{du_2},V^\epsilon_{du_1}\big] = \frac{1}{2} \, \int_s^t\int_s^{u_1} \big[v^\epsilon_{u_2},v^\epsilon_{u_1}\big]\,du_2du_1, 
$$
that is canonically associated with the space/time rescaled dynamics, equation \eqref{EqStochasticTurbulenceEq}. We put here quotation marks around "rough driver" as $V$ and $W$ only satisfy the algebraic conditions defining a rough driver, and not all of the analytic conditions since they are a priori unbounded. This is the very reason why we shall later proceed in a two step process for the analysis of the homogenisation phenomenon, by first localizing the dynamics in a ball of arbitrary radius, homogenising, and then removing the localisation.

\medskip

\subsection[\hspace{-0.6cm} {\sf Setting and result}]{Setting and result}
\label{SubsectionSetting}

Let $F$ be an almost surely continuous $\RR^d$-valued random field on $\RR^d$, defined on some probability space $(\Omega,\mcF,\PP)$. Given a measurable subset $\Lambda$ of $\RR^d$, define the $\sigma$-algebra generated by $F$ on $\Lambda$ by
$$
\cG_{\Lambda} := \sigma\big(F(x)\,; x\in \Lambda\big)\subset \mcF.
$$
We define the \textit{correlation coefficient of $F$} on two measurable subsets $\Lambda_1$ and $\Lambda_2$ of $\RR^d$ by 
$$
\alpha\big(\cG_{\Lambda_1},\cG_{\Lambda_2}\big) := \sup_{A_1\in\cG_{\Lambda_1},A_2\in\cG_{\Lambda_2}} \, \Big|\PP(A_1\cap A_2) - \PP(A_1)\,\PP(A_2)\Big|.
$$
The \textit{mixing rate of $F$} is defined as the function 
$$
\alpha(u) := \sup_{\substack{\Lambda_1,\Lambda_2\in\cB \\ \delta(\Lambda_1,\Lambda_2)\geq u}} \alpha\big(\cG_{\Lambda_1},\cG_{\Lambda_2}\big)
$$
for any non-negative $u$, and where 
$$
\delta(\Lambda_1,\Lambda_2) := \inf_{\lambda_1\in\Lambda_1,\lambda_2\in\Lambda_2}|\lambda_1-\lambda_2|.
$$ 
We make the following \textbf{\textsf{Assumptions}} on the random field $F$.  \vspace{0.15cm}

\begin{enumerate}
   \item[{\sf (i)}] \textit{The random field $F$ is centered and stationary.} \vspace{0.15cm}

   \item[{\sf (ii)}] \textit{It takes values in $\mcC^3_b(\RR^d,\RR^d)$, and }
   $$
   \EE\left[\sum_{k=0}^{3}\sup_{|x|\le 1} \big|D^k F(x)\big|^{2a_0}\right] < \infty
   $$
   \textit{for some integrability exponent} $a_0>(3\vee d)$.  \vspace{0.15cm}

   \item[{\sf (iii)}] \textit{We also have}
   $$
   \int_0^{+\infty} \alpha(u)^\kappa\,du <+\infty
   $$  
   for some exponent $\kappa\in\big(0,\frac{1}{3}\wedge \frac{1}{d}-\frac{1}{a_0}\big)$.
\end{enumerate}

\medskip

The parameters $a_0$ and $\kappa$ will be fixed throughout; we fix them once and for all. One can find in the Appendix of the work \cite{KestenPapanicolaou} of Kesten and Papanicolaou two interesting classes of examples of random fields satisfying the above assumptions, some Gaussian vector fields, and vector fields constructed from some side Poisson process.

\ssk

A last piece of notation is needed to state our main result. For any two points $x,y$ of $\RR^d$, set 
\begin{equation}
\label{EqDefnCovariance}
C(x,y) := \int_\RR \EE\Big[F(x)\otimes F(y + u{\sf v})\Big]\,du,
\end{equation}
and note that it is a function of $(y-x)$, since $F$ is stationary. This covariance function is $\mcC^{2+r}$, for any $0<u\leq 1$, under the above assumptions on $F$. One can then define a Brownian motion $V$ in the space of $\mcC^{2+r}$ vector fields on $\RR^d$, with covariance $C$, and use the results of \cite{BailleulRiedel} to define a $\mcC^{1+r}$ time-dependent random vector field $W^{(s)}_{ts}$ on $\RR^d$ by the formula
$$
W_{t,s}^{(s)}(x) = \int_s^t \int_s^{u_1} \Big[V_{\circ du_2},V_{\circ du_1}\Big](x)
$$
at each point $x$ of $\RR^d$; we use Stratonovich integration here. This can be done in such a way that the formula $\big(V,W^{(s)}+\frac{1}{2}\,V^2\big)$ defines almost surely a rough driver with regularity indices $p$ and $(2+r)$, for any $2<p<2+r<3$. Note that the integral
$$
{\sf b} := \frac{1}{2} \int_0^\infty \EE\Big[\big(D_{u{\sf v}}F\big)\,F(0) - \big(D_0F\big) \, F(u{\sf v})\Big]\,du
$$
is also well-defined as a consequence of the decay assumption {\sf (iii)} on $F$, and define a rough driver $\bV$, with regularity indices $p$ and $(2+r)$, setting
$$
{\bV}_{ts}:= \left(V_{ts} , W_{ts}^{(s)}+\frac{1}{2}\,V_{ts}^2+(t-s)\,{\sf b}\right).
$$
Denote by $\varphi$ its associated rough flow. It is proved in \cite{BailleulRiedel} that the rough flow $\varphi$ actually coincides almost surely with the stochastic flow generated by the Kunita-type stochastic differential equation 
$$
dy_t = {\sf b}\,dt + V(y_t\,; \circ dt).
$$
We read directly on this expression the generator $\mathscr{L}^{(1)}$ of the one point motion
\begin{equation*}
\label{EqGeneratorOnePoint}
\mathscr{L}^{(1)} = \overline{{\sf b}}^i\,\partial_i + \int_\RR \EE\Big[F^j(0)\,F^k(u{\sf v})\Big]\partial^2_{jk}\,du,
\end{equation*}
with 
$$
\overline{{\sf b}} := \int_0^\infty \EE\Big[\big(D_{u{\sf v}}F\big)\,F(0)\Big] du,
$$ 
and the generator $\mathscr{L}^{(2)}$ of the two point motion of the stochastic flow
\begin{equation*} 
\label{EqGeneratorTwoPoints}
\mathscr{L}^{(2)} = \mathscr{L}^{(1)}_x + \mathscr{L}^{(1)}_y + \mathscr{L}_{xy},
\end{equation*}
where $\mathscr{L}^{(1)}_x$ acts on the first variable and $\mathscr{L}^{(1)}_y$ on the second variable, while 
\begin{equation*}
\label{EqGeneratorInteraction}
\mathscr{L}_{xy}f := \int_\RR \EE\Big[F(x-y)\otimes F(x-y+u{\sf v})\Big](\partial_x,\partial_y)\,du.
\end{equation*}

\medskip

Here is how one can rephrase Kesten and Papanicolaou's homogenisation for stochastic turbulence in our flow setting.

\medskip

\begin{thm}\label{ThmHomogenisation}   {\it 
Let $K$ be any compact subset of $\RR^d$. The restriction to $K$ of $\varphi^\epsilon$ converges in law to the restriction to $K$ of $\varphi$ in $C\big([0,T]\times K\big)$.   }
\end{thm}

\medskip

We prove Theorem \ref{ThmHomogenisation} by \vspace{0.15cm}
 
\begin{enumerate}
   \item[\S \textsf{\textbf{3.2.1}}] localizing first the rough drivers ${\bfV}^\epsilon$ into a big ball of size $R$,   \vspace{0.15cm}
   
   \item[\S \textsf{\textbf{3.2.2}}] using Theorem \ref{ThmUseIt} on the localized rough drivers ${\bfV}^{\epsilon,R}$,   \vspace{0.15cm}
   
   \item[\S \textsf{\textbf{3.2.3}}] removing the localization in the end to get Theorem \ref{ThmHomogenisation}.
\end{enumerate}

\bigskip

\subsection[\hspace{-0.6cm} \textsf{Proof of homogenisation for stochastic turbulence}]{Proof of homogenisation for stochastic turbulence}
\label{SubsectionProofHomogenisation}

Recall the definition of ${\bfV}^\epsilon$ given in equation \eqref{EqDefnVep}, and let $\chi$ be a smooth real-valued function on $\RR^d$, identically equal to $1$ in the open unit ball $B(0,1)$ of $\RR^d$, with support in $B(0,2)$. Set 
$$
\chi^R(\cdot) := \chi\left(\frac{\cdot}{R}\right),
$$
and 
$$
V_{ts}^{\epsilon,R} = \chi^R V_{ts}^\epsilon =: \int_s^t v^{\epsilon,R}_u\,du  \qquad \textrm{and }\qquad  W^{\epsilon,R} := \int_s^t\int_s^{u_1} \Big[v^{\epsilon,R}_{u_2},v^{\epsilon,R}_{u_1}\Big]\,du_2du_1,
$$
and define a rough driver ${\bfV}^{\epsilon,R}$ by the formula 
$$
{\bfV}^{\epsilon,R} := \Big(V^{\epsilon,R},\,W^{\epsilon,R} + \frac{1}{2}\,\big(V^{\epsilon,R}\big)^2\Big).
$$

\bigskip

\subsubsection{Tightness of localized drivers}
\label{SubsubsectionTightness}

Our main technical ingredient for proving the tightness of the family of rough drivers ${\bV}^{\epsilon,R}$ is the following inequality due to Davydov \cite{Davydov}. If $\cA$ and $\cB$ are two sub-$\sigma$-algebras of the probability space $(\Omega,\cF)$, and if $X$ and $Y$ are two real-valued random variables that are $\cA$, resp. $\cB$, measurable, then for all integrability exponents $p_1, p_2, p_3 \geq 1$ such that $\frac{1}{p_1} + \frac{1}{p_2}+ \frac{1}{p_3} = 1$, we have
\begin{equation}
\label{EqDavydov} 
\Big| \EE[XY] - \EE[X]\,\EE[Y] \Big| \lesssim \alpha(\cA,\cB)^\frac{1}{p_1} \,\|X\|_{L^{p_2}}\,\|Y\|_{L^{p_3}}.
\end{equation}

\medskip
 
Given any $(u,x)\in\RR_+\times\RR^d$, define the set 
$$
\Lambda^\epsilon_r(x) : =x+u{\sf v} + B(0,\epsilon)
$$ 
and recall the notation $\cG_{\Lambda^\epsilon_u}$ for the $\sigma$-algebra generated by $F$ in an $\epsilon$-neighbourhood of the point $x+u{\sf v}$
 
\medskip

\begin{lem}\label{Lemma}   {\it 
\begin{itemize}
   \item[\textcolor{gray}{$\bullet$}] Let $G$ be a continuous real-valued centered random field on $\RR^d$, such that $G(r,x)\in\cG_{\Lambda^\epsilon_r(x)}$ for all positive $\epsilon$, and for which there exists a postive finite constant $m$ such that we have
$$
\EE \Bigg[\sup_{\substack{u\in[s,t],\\x \in K}} \big|G(x+u{\sf v})\big|^{2a_0}\Bigg] \le  m^{2a_0},
$$
for all $0\leq s<t\leq T$, and all compact subsets $K$ of $\RR^d$. Then we have 
$$ 
\EE\left[\left|\int_s^t G(x + u{\sf v})\,\chi_{R}(x) \, du \right|^\frac{2a_0}{a_0\kappa+1}\right]^{\frac{\kappa}{2}+\frac{1}{2a_0}} \lesssim
\chi^R(x)\, \left(\int_0^{+\infty}\alpha(u)^{\kappa}\,du \right)^\frac{1}{2} \,  \frac{m}{|{\sf v}|^{-\frac{1}{2}}}\, |t-s|^\frac{1}{2} 
$$
for all $0\leq s <t\leq T$, all $R>0$ and all $x\in \RR^d$.   \vspace{0.15cm}

   \item[\textcolor{gray}{$\bullet$}] If furthermore $H$ is another field with the same properties, and associated constant $\widetilde{m}$, then
\begin{equation}
\label{EqSecondEstimate}
\begin{split}
\EE\left[\left|\int_s^t \int_s^{u_1} G\big(x + u_1{\sf v}\big) \,  H\big(x+u_2{\sf v}\big) \, \chi_{R}(x) \, du_2 du_1 \right|^\frac{a_0}{a_0\kappa+1}\right]^{\kappa + \frac{1}{a_0}}   \\
\lesssim \left(\int_0^{+\infty}\alpha(u)^{\kappa}\,du\right)\,\chi^R(x)\, \frac{m\widetilde{m}}{|{\sf v}|}\,|t-s|. 
\end{split}
\end{equation}
\end{itemize}   }
\end{lem}
 
\medskip
 
\begin{Dem}
Set 
$$
q := \frac{a_0}{a_0\kappa+1}, \quad p_1 = 1 + \frac{1}{a_0\kappa},\quad p_2 =pa_3 = 2(a_0\kappa+1),
$$
and note that $\frac{1}{p_1} + \frac{1}{p_2} + \frac{1}{p_3} = 1$. Write first
\begin{equation*}
\begin{split}
\EE\Bigg[&\bigg|\int_s^t G(x+u{\sf v}) \chi_{R}(x)\,du \bigg|^{2q}\Bigg]^\frac{1}{2q}   \\
& \le \left(2\int_s^t \int_s^{u_1} \EE\Big[\big|G(x+u_1{\sf v}) \, G(x+u_2{\sf v}) \, \eta_R(x)^2\big|^q\Big]^\frac{1}{q} du_2du_1\right)^\frac{1}{2}    \\ 
& \lesssim \chi^R(x) \left\{ \int_s^t \int_s^r \EE\Big[\big| G(x+u_1{\sf v})\,{\bf 1}_{B(0,2R)}(x) \; G(x+u_2{\sf v})\,{\bf 1}_{B(0,2R)}(x) \big|^q\Big]^\frac{1}{q} du_2du_1 \right\}^\frac{1}{2}, 
\end{split}
\end{equation*}
and note that, for $0\leq u_2 < u_1$, we have 
$$
G(x+u_1{\sf v})\in \cG_{\Lambda^\frac{u_1-u_2}{4}_{r_1}(x)}, \quad \textrm{and} \quad  G(x+u_2{\sf v}) \textrm{ and } H(x+u_2{\sf v}) \in \cG_{\Lambda^\frac{u_1-u_2}{3}_{r_2}(x)},
$$
with
$$
d\Big(\Lambda_{u_2}^\frac{u_1-u_2}{4}(x),\Lambda_{u_1}^\frac{u_1-u_2}{4}(x)\Big) = \frac{u_1-u_2}{2}\,|{\sf v}|. 
$$
It follows from Davydov's inequality \eqref{EqDavydov} that we have 
{\small \begin{align*}
\EE\left[\left|\int_s^t G(x+u{\sf v})\,\chi^R(x) \, du\right|^\frac{2a_0}{a_0\kappa+1}\right]^{\frac{\kappa}{2} + \frac{1}{2a_0}} &\lesssim \chi^R(x)\left(\int_s^t \int_s^{u_1} \alpha\left(\frac{(u_2-u_1)|{\sf v}|}{2}\right)^\kappa m^2 \, du_2du_1\right)^{\frac{1}{2}}   \\
& \lesssim \chi^R(x) \, \left(\int_0^{+\infty} \alpha(u)^\kappa\, du\right)^{\frac{1}{2}}\, \frac{m}{|v|^{\frac{1}{2}}}\, |t-s|^{\frac{1}{2}}.
\end{align*}   }
The proof of inequaltity \eqref{EqSecondEstimate} is similar, and left to the reader.  
\end{Dem}

\medskip 

Note that the only property of $\chi$ that we used is that it has support in the ball $B(0,2)$; any derivative of $\chi$ also has this property. It follows in particular from this remark, and a change of variable, that we have 
\begin{align*}
\EE\left[\left|\eps\int_{\eps^{-2}s}^{\eps^{-2}t} (D^k_y \chi)^R \otimes D^{3-k}_{y+u{\sf v}}F \, du \right|^\frac{2a}{a\kappa+1}\right]^{\frac{\kappa}{2} + \frac{1}{2a}} &\lesssim \eps \, \Big|\eps^{-2}t - \eps^{-2}s\Big|^\frac{1}{2}\, \big|(D^k_y \chi)^R\big|  \\
& \lesssim |t-s|^{\frac{1}{2}}\,\Big|(D^k_y \chi)^R\Big|,  
\end{align*}
for all $0\leq k\leq 3$, all $y\in\RR^d$, for an implicit constant in the inequality that does not depend on $\epsilon$.

\medskip

\begin{prop}\label{theorem:tighness}   {\it 
Given a fixed positive radius $R$, the family $\big(\bV^{\eps,R}\big)_{0<\epsilon\leq 1}$ of localized rough drivers satisfies the $\epsilon$-uniform  moment bounds of \emph{\textsf{Theorem \ref{ThmUseIt}}}, with $p= 2$, any $r\in(0,1), \, a= \frac{a_0}{a_0\kappa+1}$, and for $k_1=3$.   }
\end{prop}

\medskip

\begin{DemThm}  
Starting from the representation 
{\small \begin{equation*}
\begin{split}
&\Delta^3_{\sigma}V^\eps_{t,s}(x) \\
&= \sum_{k=0}^3\binom{3}{k} \frac1{R^k}\int_{[0,1]^3}\eps\int_{\eps^{-2}s}^{\eps^{-2}t} \Big\{(D^k_{x+(u_1+u_2+u_3)\sigma}\chi)^R\otimes D^{3-k}_{x+(u_1+u_2+u_3)\sigma + u{\sf v}} \Big\}\,\sigma^{\otimes 3}\,du_1du_2du_3\,du,
\end{split}
\end{equation*}   } 
and setting 
$$
f^R := \big|\chi^R\big| + \frac1R \Big|(D\chi)^R\Big| + \frac1{R^2} \Big| (D^2\chi)^R\Big|+\frac1{R^3} \Big|(D^3\chi)^R\Big|, 
$$  
we see that 
$$
\EE\left[\Big|\Delta^3_{\sigma}V^{\eps,R}_{t,s}(x)\Big|^\frac{2a_0}{a_0\kappa+1}\right]^{\frac{\kappa}{2}+\frac{1}{2a_0}} 
\lesssim |t-s|^{\frac{1}{2}} |\sigma|^{3} \int_{[0,1]^3} f^R\Big(x+(u_1+u_2+u_3) \sigma)\Big)\,du_1du_2du_3.
$$
But since $f^R\big(x+(u_1+u_2+u_3)\sigma\big) \lesssim {\bf 1}_{B(0,R+4)}(x)$, for all $u_i\in[0,1]$ and $\sigma\in B(0,1)$, we deduce from Lemma \ref{Lemma} the inequality
\begin{equation*}
\EE\left[\Big|\Delta^3_{\sigma}V^{\eps,R}_{t,s}(x)\Big|^\frac{2a_0}{a_0\kappa+1}\right]^{\frac{\kappa}{2}+\frac{1}{2a_0}} \lesssim |t-s|^{\frac{1}{2}} |{\sf v}|^{-\frac{1}{2}}{\bf 1}_{B(0,R+4)}(x)\, |\sigma|^3.
\end{equation*}
The very same reasoning shows that the following inequalities hold 
\begin{equation*}
\begin{split}
&\EE\left[\left|\Delta^2_{\sigma} W^{\eps,R}_{t,s}(x)\right|^\frac{a_0}{a_0\kappa+1}\right]^{\kappa + \frac{1}{a_0}} \lesssim |t-s| |{\sf v}|^{-1}{\bf 1}_{B(0,R+4)}(x) \, |\sigma|^2,   \\  
&\EE\left[\left|V^{\eps,R}_{t,s}(x)\right|^\frac{2a_0}{a_0\kappa+1}\right]^{\frac{\kappa}{2}+\frac{1}{2a_0}} \lesssim |t-s|^{\frac{1}{2}} |{\sf v}|^{-\frac{1}{2}}{\bf 1}_{B(0,R+4)}(x),  \\
&\EE\left[\left| W^{\eps,R}_{t,s}(x)\right|^\frac{a_0}{a_0\kappa+1}\right]^{\kappa + \frac{1}{a_0}} \lesssim  |t-s| |{\sf v}|^{-1}{\bf 1}_{B(0,R+4)}(x);
\end{split} 
\end{equation*} 
which provide the $\epsilon$-uniform moments bounds of \textsf{Theorem \ref{ThmUseIt}} with $p= 2$, any $r\in(0,1), \, a= \frac{a_0}{a_0\kappa+1}$, and for $k_1=3$.
\end{DemThm}

\bigskip

\subsubsection{Convergence of the finite dimensional marginals of the localized drivers}
\label{SubsubsectionFiniteFimCvgce}

Fix a positive radius $R$. We know from the results of \cite{BailleulRiedel} on rough and stochastic flows that the formula
$$
Z_{ts}^{(s)} := \int_s^t \int_s^{u_1} \Big\{V_{\circ du_2} \otimes V_{\circ du_1}-V_{\circ du_1} \otimes V_{\circ du_2}\Big\},
$$
defines almost surely a $\mcC^\rho\big(\RR^d;\,L(\RR^d)\big)$-valued process. Set 
$$
{\sf c} := \frac12 \int_0^\infty \EE\Big[ F(u{\sf v})\otimes F(0) - F(0) \otimes F(u{\sf v}) \Big]\,du,
$$
and
$$
Z_{ts} := Z_{ts}^{(s)} + (t-s)\,{\sf c}.
$$
Then the formula 
\begin{equation}
\label{EqLimitRD}
\bV^R_{t,s} := \Big(V^R, (\chi^R)^2 \, W_{t,s} + \chi^R \, Z_{t,s}(\nabla \chi^R)\Big)
\end{equation}
defines almost surely a rough driver of regularity $p$ and $(2+r)$, for any $2<p<2+r<3$. We show in this section that the finite dimensional marginals of ${\bfV}^{\epsilon,R}$ converge to those of $\bV^R$ as $\epsilon$ decreases to $0$. The introduction of a notation will happen to somehow simplify our life.

\medskip

Inspecting the explicit expressions of $V^{\epsilon,R}$ and $W^{\epsilon,R}$, we see that it is sufficient to prove the finite dimensional convergence of the process
{\small $$
\left(
V^\eps_{\bullet 0};\,
DV^\eps_{\bullet 0};\,
D^2V^\eps_{\bullet 0};\,
\int_0^\bullet \big(DV^\eps_{du}\big)V^\eps_{u0};\,
\int_0^\bullet V^\eps_{du}\otimes V\epsilon_{u0};\,
\int_0^\bullet \big(D^2V^\eps_{du}\big)V^\eps_{u0};\,
\int_0^\bullet \big(DV^\eps_{du}\big)\,\big(DV^\eps_{u0}\big)
\right), 
$$   }
with values in 
$$
\RR^d\times (\RR^d)^{\otimes 2}\times (\RR^d)^{\otimes 3}\times \RR^d\times (\RR^d)^{\otimes 2}\times (\RR^d)^{\otimes 2}\times (\RR^d)^{\otimes 3},
$$
indexed by $(t,x)\in\RR_+\times\RR^d$, as $V^{\epsilon,R}$ and $W^{\epsilon,R}$ are the images of the above process by some fixed linear maps. To make that point clear, and given $\fA=\big(a^1,a^2,a^3\big)$ and $\fB=\big(b^1,b^2,b^3\big)$ in $\bigoplus_{i=1}^3(\RR^d)^{\otimes i}$, set 
$$
\fA\star \fB = \Big(a^2 b^1, a^1\otimes b^1, a^3 b^1,a^2 b^2\Big) \in \RR^d\times (\RR^d)^{\otimes 2}\times  (\RR^d)^{\otimes 2} \times (\RR^d)^{\otimes 3};
$$
remark that if $\fA\otimes \fB \in \RR^{(d+2d^2+d^3)^2}$ denotes the tensor product of $\fA$ and $\fB$, then $\fA \star \fB$ is a linear function of $\fA\otimes\fB$. If one defines the $\bigoplus_{i=1}^3(\RR^d)^{\otimes i}$-valued time-dependent fields $\frak{F}$ and $\frak{V}^\epsilon$ on $\RR^d$ by
\begin{equation*}
\begin{split}
&\frak{F}_t(x) := \Big(F(x+t{\sf v}),D_{x+t{\sf v}}F, D_{x+t{\sf v}}^2F\Big), \\
&\frak{V}^\epsilon_t(x) := \Big(V^{\epsilon}_{t0}(x), D_xV^{\epsilon}_{t0}, D^2_xV^{\epsilon}_{t0}\Big), 
\end{split}
\end{equation*}
then we have
\begin{equation*}
\begin{split}
\left(\int_0^t \big(DV^\eps_{du}\big)V^\eps_{u0}, \int_0^t V^{\eps}_{du} \otimes V^\eps_{u0} \right. &, \left. \int_0^t \big(D^2V^\eps_{du}\big) V^\eps_{u0},\int_0^t \big(DV^\eps_{du}\big) \big(DV^\eps_{u0}\big)\right) \\ 
&=  \int_0^t \mathfrak V^\eps_{du}\star \mathfrak V^\eps_u   \\
&= \epsilon^2 \, \int_{0\leq u_2\leq u_1\leq \eps^{-2}t} \frak{F}_{u_2} \star \frak{F}_{u_1}\, du_2du_1.
\end{split}
\end{equation*}
We shall use a well-paved road to prove the above finite dimensional convergence result, that can roughly be summarized as follows.

\medskip

\begin{enumerate}
   \item Decompose $\frak{V}^\epsilon$ as the sum of a martingale and small coboundary term, and use a martingale central limit theorem for dealing with the convergence of that process.   \vspace{0.1cm}
   
   \item Use a result of Kurtz and Protter on the joint convergence of pairs 
   $$
   \left(M_\bullet,\int_0^\bullet M_{u^-}dM_u\right) 
   $$
   in Skorokhod space, for good martingales, to deal with the convergence of the whole process $\Big(\frak{V}^\epsilon, \int_0^\bullet \mathfrak V^\eps_{du}\star \mathfrak V^\eps_u \Big)$. 
\end{enumerate}

\medskip

Given any positive integer $m$, any $m$-point ${\bf x} = \big(x_1,\dots,x_m\big)\in(\RR^d)^m$ in $\RR^d$, and any function $H$ from $\RR^d$ to a finite dimensional vector space, we set 
$$
H({\bf x}) := \big(H(x_1),\dots,H(x_m)\big).
$$  
We shall see in Appendix, Lemma \ref{LemAppendix}, that the following two statements are equivalent. Set $\epsilon_n = n^{-\frac{1}{2}}$.

\medskip

\begin{itemize}
   \item[\textcolor{gray}{$\bullet$}] The finite dimensional marginals of some \textit{family} of processes 
   $$
   \Big(\Big(\fV^{\eps}_\bullet,\int_0^\bullet \fV^{\eps}_{du}\star \fV^{\eps}_u\Big)(\xx)\Big)_{0<\epsilon\leq 1}
   $$ 
   converge to the corresponding finite dimensional marginals of some limit process $\big(\fV_\bullet,\frak{W}_\bullet\big)$.   \vspace{0.15cm}

   \item[\textcolor{gray}{$\bullet$}] The same convergence happens for the \textit{sequence} 
   $$
   \Big(\Big(\fV^{\epsilon_n}_\bullet,\int_0^\bullet \fV^{\epsilon_n}_{du}\star \fV^{\epsilon_n}_u\Big)(\xx)\Big)_{n\geq 1}.
   $$   \vspace{0.15cm}
\end{itemize}

We shall thus stick from now on to the study of the latter sequence. We first set the study of the sequence $\frak{V}^{\epsilon_n}$ in the setting of central limit theorems for sums of mixing, stationnary, random variables, for which martingale methods are commonly used. A positive integer $m$, and an $m$-point ${\bf x} = \big(x_1,\dots,x_m\big)\in(\RR^d)^m$, are given. Set 
$$
\lambda_{\mathbf{x}} := 2\max_{i,j} |x_i-x_j|.
$$

\medskip

\begin{lem}\label{lemma:stationarity}   {\it 
The process $\fV^1_\bullet(\xx)$ is centered, has stationary increments and is strongly mixing, with mixing coefficient $\beta$ such that 
\begin{equation}
\label{EqBoundAlphaTilde}
\beta(u) \le \alpha\big(u |{\sf v}| - \lambda_{\mathbf{x}}\big),
\end{equation}
when $u\geq \frac{\lambda_{\mathbf{x}}}{ |{\sf v}|}$.  }
\end{lem}

\medskip
 
\begin{Dem}
The spatial stationarity of the field $\frak{F}_0$({\bf x}) is inherited from the stationary character of $F$, so
$$
\Big(F({\mathbf y}+h), D_{{\mathbf y} + h}F, D^2_{{\mathbf y} + h}F\Big) \overset{\cL}{=} \Big(F({\mathbf y}), D_{{\mathbf y}}F, D^2_{{\mathbf y}}F\Big).
$$ 
It follows that 
\begin{equation*}
\fV^1_{t}(\xx)-\fV^1_{s}(\xx) = \int_{s}^{t}\fF_u(\xx) \, du \overset{\cL}{=} \int_{s}^{t}\fF_{u+h}(\xx)\,du =\fV^1_{t+h}(\xx)-\fV^1_{s+h}(\xx),
\end{equation*}
so the process $\fV^1_\bullet(\xx)$ has stationary time-increments. As we also see on the first equality that the random variable $\big(\fV^1_{t}(\xx)-\fV^1_{s}(\xx)\big)$ is $\cG_{[s{\sf v},t{\sf v}]+B(0,\lambda_{\mathbf{x}})}$-measurable, bound \eqref{EqBoundAlphaTilde} follows from the inequality
$$
\delta\Big([s_1{\sf v} , t_1 {\sf v}] + B(0,\lambda_{\mathbf{x}})\, , \,[s_2{\sf v} , t_2 {\sf v}] + B(0,\lambda_{\mathbf{x}})\Big) \ge (s_2-t_1)|{\sf v}| - \lambda_{\mathbf{x}}, 
$$
which holds for all $s_1\le t_1 \le s_2 \le t_2$, with $(s_2-t_1) |{\sf v}|>\lambda_{\mathbf{x}}$.
\end{Dem}
 
\medskip

Define a stationary sequence of fields setting 
$$
\mathfrak X_k := \int_k^{k+1} \mathfrak F_u \, du.
$$ 
We shall analyse the asymptotic behaviour of $\left(\fV^{\epsilon_n}_t,\int_0^t\fV^{\epsilon_n}_{\dd r} \star \fV^{\epsilon_n}_r\right)(\xx)$ by first writing it in terms of the $\mathfrak X_k$, in the next lemma, and then by using a "martingale plus remainder" decomposition -- see Lemma \ref{lemma:bound_theta}. This will then put us in a position to use a well-known result of Kurtz and Protter about the convergence in law of pairs of the form $\big(M_\bullet,\int_\bullet MdM\big)$, for good martingales $M$. 

\medskip

\begin{lem}\label{lemma:discretization2}   {\it
Given any time $t\in [0,T]$, any positive integer $m$, and any $m$-point $\xx$, we have
\begin{multline*}
\left(\fV^{\epsilon_n}_t,\int_0^t\fV^{\epsilon_n}_{\dd r} \star \fV^{\epsilon_n}_r\right)(\xx) = 
\left(\frac1{\sqrt n} \sum_{k=0}^{[nt]-1}\fX_k,\frac1n\sum_{k=0}^{[nt]-1} \sum_{j=0}^k \mathfrak X_k \star \mathfrak X_j\right)(\xx) \\
 + \left(0,t\int_{0\leq u_2\leq u_1\leq 1} \EE\big[\mathfrak F_{u_1}\star \mathfrak F_{u_2}\big]\,du_2du_1\right)(\xx) + R^n_t(\xx)
\end{multline*} 
with a remainder $R^n_t(\xx)$ that converges to $0$ almost surely.   }
\end{lem}

\medskip

So Slutsky's theorem brings back the study of the finite dimensional convergence of the process in the left hand side of the above equality to the study of the finite dimensional convergence of the first term in the right hand side. 

\medskip

\begin{Dem}
Let first pick an $R$ bigger than all the $|x_i|$, and note that $\frak{V}^\epsilon_\bullet = \frak{V}^{\epsilon,R}_\bullet =: \Big(V^{\epsilon,R},DV^{\epsilon,R},D^2V^{\epsilon,R}\Big)$, and $\int_0^\bullet\frak{V}^\epsilon_{du}\star\frak{V}^\epsilon_u\,du = \int_0^\bullet\frak{V}^{\epsilon,R}_{du}\star\frak{V}^{\epsilon,R}_u\,du$, on $B(0,R)$. As it follows from the tightness result stated in theorem \ref{theorem:tighness} that there is almost surely an exponent $b<\frac{1}{2}$ such that the processes $\frak{V}^{\epsilon,R}_\bullet = \frak{V}^\epsilon_\bullet$ have finite $b$-H\"older norm uniformly in $\epsilon$, and since 
$$
\frak{V}^{\epsilon_n}_{nt}(\xx) - \frac{1}{\sqrt{n}}\,\sum_{k=1}^{[nt]}\frak{X}_k(\xx) = \Big(\frak{V}^{\epsilon_n}_{nt}(\xx)-\frak{V}^{\epsilon_n}_{[nt]}(\xx)\Big)
$$
the first component of $R_n(\xx)$ converges indeed to $0$ almost surely. To analyse its second component, write 
\begin{align*}
\int_0^t\fV^{\epsilon_n}_{\dd r} \star \fV^{\epsilon_n}_r(\xx) &= \frac{1}{n} \, \int_{0\leq u_2\leq u_1\leq nt} \mathfrak F_{u_1} \star \mathfrak F_{u_2}\,du_2du_1  \\
&= \frac{1}{n} \, \int \Big\{ {\bf 1}_{0\leq u_2\leq u_1\leq [nt]} + {\bf 1}_{0\leq u_1\leq [nt]}{\bf 1}_{[u_1]\leq u_2\leq u_1} \\
&\quad+ {\bf 1}_{[nt]\leq u_1\leq nt}{\bf 1}_{0\leq u_2\leq [u_1]} + {\bf 1}_{[nt]\leq u_1\leq nt}{\bf 1}_{[u_1]\leq u_2\leq u_1}\Big\}\,\mathfrak F_{u_1} \star \mathfrak F_{u_2}\,du_2du_1.
\end{align*}
\begin{itemize}
   \item[\textcolor{gray}{$\bullet$}] The first term is equal to 
   $$
   \frac1n\sum_{k=1}^{[nt]-1} \sum_{j=0}^{k-1} \mathfrak X_k \star \mathfrak X_j = \frac1n\sum_{k=0}^{[nt]-1} \sum_{j=0}^{k-1} \mathfrak X_k \star \mathfrak X_j + 0\Big(\frac{1}{n}\Big). \vspace{0.1cm}
   $$
   
   \item[\textcolor{gray}{$\bullet$}] Writing the second term as
   $$
   \frac{1}{n}\,\sum_{k=0}^{[nt]-1} \left(\int_{k\leq u_2\leq u_1\leq k+1}\mathfrak F_{u_1}\star \mathfrak F_{u_2}(\xx)\,du_2du_1\right),  
   $$
   it appears as the ergodic sum of the stationary mixing sequence given by the integral term. The ergodic theorem gives its asymptotic behaviour.   \vspace{0.1cm}

   \item[\textcolor{gray}{$\bullet$}] The third term decouples and writes
   $$
   \left(\int_{[nt]}^{nt} \fF_u\,du\right) \star\left(\frac{1}{n} \int_0^{[nt]}\mathfrak F_u\,du\right) (\xx) = O\Big(\frac{1}{n}\Big)\,\Big(t\,\EE\big[\frak{F}_1\big]+o_n(1)\Big) = O\Big(\frac{1}{n}\Big)\,o_n(1),    
   $$
   by Lemma \ref{lemma:stationarity} and the ergodic theorem.   \vspace{0.1cm}
      
   \item[\textcolor{gray}{$\bullet$}] The fourth term is almost surely of order $\frac{1}{n}$.
\end{itemize} 
\end{Dem}
 
\medskip
 
To set the scene of Gordin's martingale decomposition fof stationary sequences, define  
$$
\cF_k(\xx) := \sigma\Big(\fX_j(\xx), j\leq k\Big)
$$
and
$$
\theta_k(\xx) := \sum_{j\geq 0} \EE\Big[\mathfrak X_{k+j}(\xx)\Big | \cF_k\Big] = \int_k^\infty \EE\Big[\mathfrak F_u(\xx) \Big|\cF_k\Big] \,du.
$$
The fact that the sequence $\big(\fX_k(\xx)\big)_{k\geq 0}$ is stationary and mixing, with mixing coefficient $\beta(k)$, ensures that $\theta_k(\xx)\in L^2$, so it is in particular almost surely finite. Note the relation 
$$
\mathfrak X_k(\xx) = \theta_k(\xx)-\theta_{k+1}(\xx) +\Big(\theta_{k+1}(\xx)-\EE[\theta_{k+1}(\xx)|\cF_k]\Big).
$$
Denote by  $\big(\mathfrak M_k(\xx)\big)_{k\geq 0}$ the $L^2$-martingale with increments $\theta_{k+1}(\xx)-\EE\big[\theta_{k+1}(\xx) \big| \cF_k\big]$; so
$$
\mathfrak M_k(\xx) - \fM_0(\xx) = \sum_{j=0}^{k-1}  \fX_j(\xx) + \theta_k(\xx)-\theta_0(\xx)= \int_0^\infty \Big\{\EE\big[\mathfrak F_u(\xx) \big|\cF_k\big] - \EE\big[\mathfrak F_u(\xx) \big|\cF_0\big]\Big\} \, du. 
$$
We define a  pure jump c\`adl\`ag martingale by the formula
$$
\mathfrak M^n_t(\xx) = \frac1{\sqrt n} \mathfrak M_{[nt]}(\xx);
$$
it satisfies the relation 
$$
\int_0^t \fM^n_{du}\star\fM^n_u = \frac{1}{n} \sum_{k=0}^{[nt]-1}(\fM_{k+1}-\fM_k)\star\fM_k,
$$
with an It\^o integral used in the left hand side.

\medskip

\begin{lem}\label{lemma:bound_theta}   {\it
We have
\begin{align*}
\left(\frac1{\sqrt n}\sum_{k=0}^{[nt]-1}\fX_k,\frac1n\sum_{k=0}^{[nt]-1}\sum_{j=0}^{k}\fX_k\star\fX_j\right)&(\xx) = \left( \fM^n_t,\int_0^t \fM^n_{\dd r}\star \fM^n_r \right)(\xx)   \\
& + \left(0, t\, \int_{[1,\infty)\times [0,1]} \EE\big[\fF_{u_1}\star\fF_{u_2}(\xx)\big]\,du_2du_1\right) + \overline{R}^n_t, 
\end{align*}
for a remainder $\overline{R}^n_t$ that converges in probability to $0$.   }
\end{lem}
 
\medskip
 
\begin{Dem}
{\bf (1)} We start giving a uniform bound on $\theta_k(\xx)$. Define for that purpose an exponent $q$ by the relation $\kappa + \frac{a\kappa+1}{2a} + \frac{1}{q} = 1$, and let $Q\in L^q\cap \cF_k(\xx)$ have unit $L^q$-norm. We have the $k$-uniform bound
\begin{align*}
\Big|\EE\big[Q\theta_k(\xx)\big]\Big|
&\lesssim 
\int_k^\infty \Big|\EE\big[Q\fF_u(\xx)\big]\Big| \, du \\
&\lesssim
\int_k^\infty \alpha\Big(\big| u-(k+1)\big|{\sf v}|\big|\Big)^\kappa \,du \\
&\lesssim_\alpha \, 1, 
\end{align*}
so $\theta_k(\xx)$ have a finite $L^\frac{2a_0}{3a_0\kappa+1}$-norm, uniformly bounded as $k$ varies.

\ssk

{\bf (2)} \textcolor{gray}{$\bullet$} The case of 
\begin{align*}
\frac1{\sqrt n} \sum_{k=0}^{[nt]-1} \fX_k(\xx) -\fM_t^n(\xx) 
& =
\frac{1}{\sqrt n} \big( \theta_0(\xx)- \theta_{[nt]}(\xx)\big) + \frac{1}{\sqrt n} \fM_0(\xx).
\end{align*}
is trivially dealt with using the above $k$-uniform boundon $\theta_k(\xx)$.

\textcolor{gray}{$\bullet$} For the second component, start with the decomposition
{\small \begin{equation}
\label{EqProofMartingale}
\begin{split}
 \frac1n &\sum_{k=1}^{[nt]-1} \sum_{j=0}^k \fX_k \star \fX_j(\xx)
 -\frac1n \sum_{k=0}^{[nt]-1} (\fM_{k+1}-\fM_k)\star \fM_k(\xx)\\
  =&
  \frac1n\sum_{k=1}^{[nt]-1} \fX_{k}\star(\fM_{k} - \theta_{k} + \theta_0-\fM_0)(\xx) - \frac1n\sum_{k=0}^{[nt]-1} (\fX_{k} + \theta_{k+1}-\theta_k)\star \fM_k(\xx) \\
  =& - \frac1n\sum_{k=0}^{[nt]-1} \fX_{k}\star \theta_k(\xx) - \frac1n\sum_{k=0}^{[nt]-1}(\theta_{k+1}-\theta_{k})\star\fM_k(\xx) + \Big(\frac1n\sum_{k=1}^{[nt]-1} \fX_{k}\Big)\star(\theta_0-\fM_0)(\xx).
\end{split}
\end{equation}   } 
Since the centered sequence $\big(\fX_k(\xx)\big)_{k\geq 0}$ is stationary and mixing,the last term above converges to $0$, by the ergodic theorem. The sequence $\big(\big(\fX_k\star\theta_k\big)(\xx)\big)_{k\geq 0}$ is also stationary and ergodic, so $- \frac1n\sum_{k=0}^{[nt]-1} \fX_{k}\star \theta_k(\xx)$ converges almost surely to $-\EE\big[\fX_0\star \theta_0 (\xx)\big]$. To analyse the second term in the right hand side of equation \eqref{EqProofMartingale}, write it as 
\begin{equation*}
- \frac1n\sum_{k=0}^{[nt]-1}\sum_{j=0}^{k-1}(\theta_{k+1}-\theta_{k})\star(\fM_{j+1}-\fM_{j})(\xx) + \frac1n(\theta_{[nt]}-\theta_0)\star\fM_0(\xx).
\end{equation*}
The second term converges to zero in probability, by point \textbf{(1)}. Remark that 
\begin{align*}
- \frac1n\sum_{k=0}^{[nt]-1}\sum_{j=0}^{k-1}\big(\theta_{k+1}-\theta_{k}\big)\star\big(\fM_{j+1}-\fM_{j}\big)(\xx) 
&=  \frac1n\sum_{j=0}^{[nt]-2}\theta_{j+1}\star \big(\fM_{j+1}-\fM_{j}\big)(\xx)  \\
&- \theta_{[nt]-1}\star\left(\frac1n\sum_{k=0}^{[nt]-2}\fM_{k+1}-\fM_{k}\right)(\xx). 
\end{align*}
Here again, thanks to the ergodic theorem the second term in the right hand side converges to zero in probability. Furthermore, by construction, the sequence $\big(\theta_{j+1}\star(\fM_{j+1}-\fM_{j})(\xx)\big)_{j\geq 0}$ is stationary and ergodic, so the first term of the right hand side converges almost surely to $\EE[\theta_1\star (\fM_1-\fM_0)(\xx)]$. All these elementary remarks together prove that
\begin{align*}
\frac1n \sum_{k=1}^{[nt]-1} \sum_{j=0}^k \fX_k \star \fX_j(\xx)
&-\frac1n \sum_{k=0}^{[nt]-1} \big(\fM_{k+1}-\fM_k\big)\star \fM_k(\xx)
\end{align*}
converges in probability to
$$
-\EE\Big[\fX_0\star \theta_0 (\xx)\Big] + \EE\Big[\theta_1\star (\fM_1-\fM_0)(\xx)\Big].
$$ 
 
{\bf (3)} In order to prove the lemma, it remains to find a good expression for the limit. For all $j\ge 1$, we have
 $$
\fX_j\star \fX_0(\xx) - \big(\fM_{j+1}-\fM_j\big) \star \big(\fM_{1}-\fM_0\big) = \fX_j \star \big(\theta_0-\theta_1\big)(\xx) + \big(\theta_j - \theta_{j+1}\big)\star \big(\fM_1 - \fM_0\big)(\xx),
$$
with 
$$
\EE \Big [ \big(\fM_{j+1}-\fM_j\big) \star \big(\fM_{1}-\fM_0\big)\Big] = 0,
$$ 
since $\fM(\xx)$ is an $\big(\cF_k(\xx)\big)_{k\geq 0}$-martingale. One can then use the fact that 
$$
\EE\Big[\fX_j \star \big(\theta_0 - \theta_1\big) (\xx)\Big] = \EE\Big[\fX_{j+\ell}\star \big(\theta_\ell-\theta_{\ell+1}\big) (\xx)\Big]
$$
for all $j,\ell\geq 1$, to see that 
{\small \begin{align*}
\sum_{j=1}^N \EE\Big[\fX_j \star \fX_0 (\xx)\Big] 
&= \sum_{j=1}^N \EE\Big[\fX_j \star \fX_0 (\xx)\Big] - \EE \Big[ \big(\fM_{j+1}-\fM_j\big) \star \big(\fM_{1}-\fM_0\big)\Big]\\
&= \sum_{j=1}^N \EE \Big[\fX_j \star \big(\theta_0-\theta_1\big)(\xx) + \big(\theta_j - \theta_{j+1}\big)\star \big(\fM_1 - \fM_0\big)(\xx) \Big]\\
&= \sum_{j=1}^N \EE\Big[\fX_N \star \big(\theta_{N-j} - \theta_{N-j+1}\big)(\xx) + \big(\theta_j - \theta_{j+1}\big) \star \big(\fM_1 - \fM_0\big)(\xx)\Big]\\
&= -\EE\Big[\fX_N \star \theta_{N}(\xx)\Big] + \EE\Big[\theta_1 \star \big(\fM_1 - \fM_0\big) (\xx)\Big] + \textrm{R}^N\\
&= -\EE\Big[\fX_0 \star \theta_{0}(\xx)\Big] + \EE\Big[\theta_1 \star \big(\fM_1 - \fM_0\big) (\xx)\Big] + \textrm{R}^N,
\end{align*}   }
where 
$$
\textrm{R}^N := \EE\Big[\fX_N\star \theta_0-\theta_{N+1}\star \big(\fM_1 - \fM_0\big)(\xx)\Big]
$$ 
converges to zero as $N$ goes to infinity, thanks to the mixing properties of $\fX$ and $\theta$. It follows that
\begin{align*}
\frac1n \sum_{k=1}^{[nt]-1} \sum_{j=0}^k \fX_k \star \fX_j(\xx) - \frac1n \sum_{k=0}^{[nt]-1} \big(\fM_{k+1}-\fM_k\big)\star \fM_k(\xx)
\end{align*}
converges in probability to
\begin{align*}
\sum_{j=1}^\infty \EE\big[\fX_j \star \fX_0 (\xx)\big] &= \sum_{j=1}^\infty \int_j^{j+1}\int_0^1 \EE\Big[\cF_{u_1} \star \cF_{u_2} (\xx)\Big] \, du_2du_1 \\
&= \int_1^\infty\int_0^1 \EE\Big[\cF_{u_1} \star \cF_{u_2}(\xx)\Big] \, du_2du_1.
\end{align*} 
\end{Dem}

\medskip
 
We are now ready to prove the finite dimensional convergence of the sequence of processes $\Big(\frak{V}^{\epsilon_n}_\bullet,\int_0^\bullet \frak{V}^{\epsilon_n}_{du}\star \frak{V}^{\epsilon_n}_u\Big)$ to the process $\Big(\frak{V}_\bullet,\int_0^\bullet \frak{V}_{du}\star \frak{V}_u\Big)$, where $\frak{V}$ is a Brownian motion on the space $\Big(\RR^d\oplus(\RR^d)^{\otimes 2}\oplus(\RR^d)^{\otimes 3}\Big)^m$, with covariance ${\sf C}(\xx)$ given, for all $\lambda,\mu\in\Big(\RR^d\oplus(\RR^d)^{\otimes 2}\oplus(\RR^d)^{\otimes 3}\Big)^m$, by the formula
\begin{equation*}
\begin{split}
\big({\sf C}(\xx) \lambda\big)\mu &= \sum_{k\in\ZZ}\EE\Big[\big(\fX_0(\xx)\cdot\lambda\big)\big(\fX_k(\xx)\cdot\mu\big)\Big]   \\
&= \int_0^1\int_\RR \EE\Big[\big(\fF_{u_1}(\xx)\cdot\lambda\big) \big(\fF_{u_2}(\xx)\cdot\mu\big) \Big]\,du_1du_2   \\
&= \int_\RR \EE\Big[\big(\fF_0(\xx)\cdot\lambda\big) \big(\fF_{r}(\xx)\cdot\mu\big) \Big]\,du.  
\end{split}
\end{equation*}
(We used the time-stationarity of $\frak{F}$ in the last line.) We shall use for our purposes a useful result proved by Kurtz and Protter in \cite{KurtzProtter}, that says that if $(M^n)_{n\geq 1}$ is a sequence of vector-valued martingales, with $\EE\big[M^n_1\big]=0$ for all $n\geq 0$, and if $\big(M_t^n\big)_{n\geq 0}$ is bounded in $L^2$ for each time $t\in [0,T]$, then the convergence in law of $M^n_\bullet$ to $M_\bullet$ in Skorokhod space implies the convergence in law of the pair $\big(M^n_\bullet, \int_0^\bullet M^n_{u^-} \, dM^n_u \big)$ to $\big(M_\bullet, \int_0^\bullet M_{u^-} \, dM_u\big)$, in the Skorokhod space -- the integrals are understood as It\^o integrals. 

\ssk

Indeed, since Lemma \ref{lemma:stationarity} ensures that we can use for the sequence of processes $\frac{1}{\sqrt{n}}\sum_{k=1}^{[n\bullet]}\frak{X}_k$ well-known invariance principles, Lemma \ref{lemma:bound_theta} shows that the process $\frak{M}^n_\bullet(\xx)$ converges to the process $\frak{V}_\bullet$ in Skorokhod space. As moreover the sequence $\big(\theta_k(\xx)\big)_{k\geq 0}$ is bounded in $L^\frac{2a_0}{3a_0\kappa+1}$, and $\kappa<\frac{1}{3}-\frac{1}{a_0}$, that sequence is also bounded in $L^2$. So we have 
\begin{equation*}
\begin{split}
\EE\Big[\Big|\frak{M}^n_t(\xx)\Big|^2\Big] &\lesssim \frac{1}{n}\,\EE\left[\left| \sum_{k=0}^{[nt]-1}\frak{X}_k(\xx)\right|\right] + \frac{1}{n}\,\EE\big[|\theta_0|^2\big] + \frac{1}{n}\,\EE\big[|\theta_{nt}|^2\big] \\
&\lesssim t + O\Big(\frac{1}{n}\Big),
\end{split}
\end{equation*} 
which implis that the random variables $\frak{M}^n_t$ are uniformly bounded in $L^2$, for all $t\in [0,T]$ and $n\geq 1$. This fact finally puts us in a position to use the above mentioned result of Kurtz and Protter. 
 
\bigskip

Putting all pieces together, we have proved that the finite dimensional laws of the family of processes 
$$
\Big(\fV^{\eps,R}_\bullet(\xx),\int_0^\bullet \fV^{\eps,R}_{du}\star \fV^{\eps,R}_u (\xx) \Big)
$$
converge weakly to the finite dimensional laws of the process 
$$
\Big(\fV_\bullet(\xx), \int_0^\bullet \fV_{du}\star \fV_u (\xx) + \bullet \mathfrak b(\xx)\Big),
$$
where
\begin{equation*}
\begin{split}
\mathfrak b(\xx) &:= \int_0^1\int_0^{u_1}\EE\Big[\fF_{u_1}\star \fF_{u_2}(\xx)\Big]\,du_2du_1 + \int_1^\infty \int_0^1\EE\Big[\fF_{u_1}\star \fF_{u_2}(\xx)\Big]\,du_2du_1  \\ 
						    &= \int_0^\infty \EE\Big[\frak{F}_u\star\frak{F}_0(\xx)\Big]\,du = \int_0^\infty \EE\Big[\frak{F}_u\star\frak{F}_0\Big]\,du =: \frak{b}
\end{split}
\end{equation*}
is independent of $\xx$, by stationarity of $\frak{F}$. Explicitly, one can write
$$
\mathfrak b = \mathfrak b(\xx) = \big(\underbrace{{\sf b}^1,{\sf b}^2,{\sf b}^3,{\sf b}^4,\cdots,{\sf b}^1,{\sf b}^2,{\sf b}^3,{\sf b}^4}_{m\ times}\big),
$$
with 
\begin{align*}
{\sf b}^1 &:= \int_0^\infty \EE\Big[\big(D_{u{\sf v}}F\big)\,F(0)\Big] \, du, \qquad {\sf b}^2 := \int_0^\infty \EE\Big[F(u{\sf v}) \otimes F(0)\Big] \, du  \\
{\sf b}^3 &:= \int_0^\infty \EE\Big[\big(D^2_{u{\sf v}}F\big)\, F(0)\Big] \, du,\qquad {\sf b}^4 := \int_0^\infty \EE\Big[D_{u{\sf v}}F\,D_0F\Big] \,du. \\
 \end{align*}
Write $V$ for the first component of $\frak{V}$. If one recalls now that ${\bV}^{\epsilon,R}$ is obtained from $\Big(\fV^{\eps,R}_\bullet(\xx),\int_0^\bullet \fV^{\eps,R}_{du}\star \fV^{\eps,R}_u (\xx) \Big)$ by a fixed linear map, it follows that the finite dimensional laws of ${\bV}^{\epsilon,R}$ converge to the finite dimensional laws of the rough driver
\begin{align*}
\widehat{\bV}^R_{ts}(x) := \bigg(&\chi^R(x) V_{t,s}(x)\, , \, \chi^R(x)^2 \left\{\int_s^t \big(D_xV_{du}\big)V_{us}(x) - \frac{1}{2}\,\big(D_xV_{ts}\big)V_{ts}(x) + (t-s)\,{\sf b}^1 \right\}   \\
&+ \chi^R(x) \left\{\int_s^t V_{du}(x)\otimes V_{us}(x) -\frac12 V_{t,s}(x)\otimes V_{t,s}(x) + (t-s)\, {\sf b}^2\right\}\nabla\chi^R(x)\bigg).
\end{align*}
In its Stratonovich form, and with the notations introduced before theorem \ref{theorem:convergence_KP}, this gives
\begin{align*}
\widetilde \bV^R_{ts}(x) = \bigg(\chi^R(x) V_{t,s}(x)\, , \, \chi^R(x)^2 &\Big\{W^{(s)}_{ts}(x) + (t-s)\big({\sf b}^1 - \partial_1C(0,0)\big) \Big\}  \\
&+ \chi^R(x) \Big\{Z^{(s)}_{ts}(x) + (t-s) \big({\sf b}^2 - C(0,0)\big)\Big\}\nabla \chi^R(x)\bigg),
\end{align*}
where $\partial_1C(0,0) = \partial_x C(x,y)_{|x=y=0}$. Since
$$
{\sf b}^1 - \partial_1C(0,0) = \frac12 \int_0^\infty \EE\big[\big(D_{u{\sf v}}F\big) F(0) - \big(D_0F\big) F(u{\sf v})\big]\,du = {\sf b},  
$$
and
$$
{\sf b}^2 - C(0,0) = \frac12 \int_0^\infty \EE\big[F(u{\sf v})\otimes F(0) - D_0F\otimes F(u{\sf v})\big]\,du = {\sf c},
$$
we finally see that 
$$
\widehat\bV^R = \bV^R.
$$ 
This fact finishes the proof of the convergence of the finite dimensional laws of ${\bV}^{\epsilon,R}$ to those of ${\bV}^R$.

\bigskip

\subsubsection{End of proof of Theorem \ref{ThmHomogenisation}}
\label{SubsubsectionEndOfProof}

The results of sections \ref{SubsubsectionTightness} and \ref{SubsubsectionFiniteFimCvgce} together put us in a position to use Theorem \ref{ThmUseIt} and show that the family of rough drivers ${\bfV}^{\epsilon,R}$, with regularity indices $p$ and $(2+r)$, converges as a $(p',2+r')$-rough driver to the rough driver ${\bfV}^R$ with regularity indices $p'$ and $(2+r')$, introduced in \eqref{EqLimitRD}. One finally uses the following elementary fact to remove the localisation.

\medskip

\begin{prop}\label{theorem:loc}   {\it 
Assume that the quantities
$$
\EE\Big[\sup_{0\leq s\leq t\leq T}\,\sup_{x\in K}\,\big| \varphi^{\epsilon,R}(x)\big|\Big]
$$
are uniformly bounded above by a constant independent of $R$, for each compact subset $K$ of $\RR^d$. Then, the restriction $\varphi^\epsilon_{| K}$ of $\varphi^\epsilon$ to $K$ converges in law in $C(D_T\times K)$ to the restriction $\varphi_{| K}$ of $\varphi$ to $K$.   }
\end{prop}

\medskip

\begin{Dem}
Given a compact subset $K$ of $\RR^d$, set
$$
M := \EE\Big[\sup_{0\leq s\leq t\leq T}\,\sup_{x\in K}\,\big| \varphi^{\epsilon,R}(x)\big|\Big] < \infty.
$$
Given any closed set $F$ of $C(\Delta_T\times K)$, we have
\begin{equation*}
\begin{split}
\PP\Big(\varphi^\epsilon_{| K} \in F\Big) &\leq \PP\Big(\varphi^{\epsilon,R}_{| K} \in F\Big) + \PP\Big(\varphi^\epsilon(K)\cap B(0,R)^c\neq\emptyset\Big)  \\
&\leq \PP\Big(\varphi^{\epsilon,R}_{| K} \in F\Big) + \PP\Big(\varphi^{\epsilon,R}(K)\cap B(0,R)^c\neq\emptyset\Big),   
\end{split}
\end{equation*}
with 
$$
\underset{\varepsilon}{\limsup}\, \PP\Big(\varphi^{\epsilon,R}_{| K} \in F\Big) \leq \PP\Big(\varphi^R_{| K} \in F\Big)
$$
by the convergence assumption on the rough driver ${\bfV}^{\epsilon,R}$ and the continuity of the It\^o map, while the second term can be bounded above by 
$\frac{M}{R}$. The conclusion follows by letting $R$ tend to $\infty$. 
\end{Dem}

\bigskip
\bigskip

\appendix

\section[\hspace{0.6cm} {\sf Compactness results}]{Compactness results for rough drivers}
\label{SectionAppendixCompactness}

We prove in this section a Lamperti-type compactness criterion for random rough drivers that implis Theorem \ref{ThmUseIt} in a straightforward way. We shall use for that purpose an elementary result on Besov spaces which we recall first.

\medskip

Given $f\in L^\infty(\RR^d,\RR^d)$  and $\sigma$ in the unit ball of $\RR^d$, we define inductively a sequence $\Delta_\sigma^m$ of operators on $L^\infty(\RR^d,\RR^d)$ setting
$$
\big(\Delta_\sigma f\big)(\cdot) = f(\cdot + \sigma) - f(x) \quad \mathrm{and} \quad \Delta^{m+1}_\sigma f = \Delta_\sigma (\Delta^m_\sigma f).
$$
Given positive parameters $a,b\leq \infty$, and two exponents $0<\alpha\leq m$, the Besov space $B^\alpha_{ab}(\RR^d) =: B^\alpha_{ab}$ is defined as 
$$
\left\{ f\in L^p(\RR^d) \, : \,  \|f\|_{B^\alpha_{ab}} := \|f\|_{L^a(\RR^d)} + \left( \int_{B(0,1)} |\sigma|^{-b \alpha }\|\Delta^m_\sigma f\|^b_{L^a(\RR^d)} \frac{\dd \sigma}{|\sigma|^d}\right)^\frac{1}{b}< +\infty \right\}.
$$
Two different choices of constants $m (\geq \alpha)$ define the same space, with equivalent norms; so we do not keep track of that parameter in the notation for the space. These spaces provide refinements of the H\"older spaces, in so far as $B_{\infty,\infty}^\alpha = \mcC^\alpha$, for non-integer $\alpha$'s. The most useful property of this scale of spaces will be for us Besov's embedding properties, according to which, if one is given $1\leq p_1 < p_2\leq \infty$ and $\alpha>0$, then $B^{\alpha + d\big(\frac{1}{p_1} - \frac{1}{p_2}\big)}_{p_1,p_1}$ is continuously embedded into $B^\alpha_{p_2,p_2}$. The following elementary continuity result was also used above.

\medskip

\begin{prop}\label{prop:product}   {\it 
Let $0<\alpha_1\leq \alpha_2$ the multiplication is a continuous bilinear operator from $\cC^{\alpha_1}\times \cC^{\alpha_2}$ to $\cC^{\alpha_1}$.   }
\end{prop} 

\medskip

\noindent From a probabilistic point of view, the interest of working with Besov spaces comes from the fact that it is usually hard to get estimates on the expectation of some supremum, while making computions on integral quantities is usually much easier, as the proof of the next proposition will make it clear. We use in this statement the notations
$$
\nu_V(d\sigma) := \mu^{-((2+r)a + d)}\,d\sigma, \qquad \nu_W(d\sigma) := \mu^{-((r+1)a + d)}\,d\sigma,
$$
for two measures on the unit ball $B(0,1)$ of $\RR^d$, absolutely continuous with respect to Lebesgue measure $d\sigma$, and for a range of parameters $(a,r)$ specified in the statement. Recall $D_T$ stands for the $2$-dimensional simplex $\big\{(s,t)\in[0,T]^2\,;\,s\leq t\big\}$.

\medskip

\begin{prop}\label{theorem:kolmo}   {\it 
Assume we are given a family $\big(V_{ts},W_{ts}\big)_{0\leq s\leq t\leq T}$ of random vector fields on $\RR^d$, with $V$ almost surely additive as a function of time, of class $\mcC^2$, and with $W$ satisfying almost surely the identity
\begin{equation}
\label{eq:alg2} 
W_{ts} = W_{tu} + W_{us} + \frac12\big[V_{us},V_{tu}\big],
\end{equation}
for every $0\leq s\leq t\leq T$. Let $a,p$ and $r$ be positive parameters, with $a\geq 1$ and 
$$
0 < \frac{1}{\frac{1}{p} - \frac{1}{2a}} - 2 < r - \frac{d}{a} < 1.
$$ 
Assume that there exists two non-negative functions $C_0^V\in L^a(\RR^d)$ and $C_0^W\in L^a(\RR^d)$ such that we have 
\begin{equation}\label{eq:kolmo1}
\left\|\frac{V_{t,s}(y)}{|t-s|^\frac{1}{p}}\right\|_{L^{2a}} \leq C_0^V(y), \qquad\textrm{and}\qquad \left\|\frac{W_{t,s}(y)}{|t-s|^\frac{2}{p}}\right\|_{L^a} \leq C_0^W(y)
\end{equation}
for all $(y,(s,t))\in \RR^d\times D_T$. Assume also that there exists an integer $k_1\geq 3$, and two functions 
$$
C_1^V \in L^a \Big(\big(B(0,1),\mu_V\big) ;  L^a(\RR^d)\Big),
$$ 
and 
$$
C_1^W\in L^a \Big(\big(B(0,1),\mu_W\big) ;  L^a(\RR^d)\Big),
$$ 
such that we have 
\begin{equation}\label{eq:kolmog}
\begin{aligned}
\left\|\frac{\Delta_\sigma^{k_1} V_{t,s}(y)}{|t-s|^\frac{1}{p}}\right\|_{L^{2a}} \leq  C_1^V(\sigma\,;y), \qquad\textrm{and}\qquad \left\|\frac{\Delta_\sigma^{k_1-1} W_{t,s}(y)}{|t-s|^\frac{2}{p}}\right\|_{L^a} \leq  C_1^W(\sigma\,;y),
\end{aligned}
\end{equation} 
for all $(\sigma,y,(s,t))\in B(0,1)\times \RR^d \times D_T$. Then, for any $(p',2+r')$ with $p'< 2 + r'$ and 
$$
r'<r-\frac{d}{a},\qquad \textrm{ and }\qquad  \frac{1}{3}<\frac{1}{p'}<\frac{1}{p}-\frac{1}{2a},
$$
there exists a modification $\widetilde{\bV}$ of ${\bV} := \big(V,\frac{1}{2}V^2+W\big)$ that is almost surely a rough driver with regularity indices $p'$ and $(2+r')$, for which 
$$
\EE\left[\big\|\widetilde{\bV}\big\|_{(p',2+r')}^{2p}\right] \lesssim \big\|C_0^V\big\|+ \big\|C_0^W\big\| + \big\|C_1^V\big\| + \big\|C_1^W\big\|,
$$
with each norm taken in its natural space.   }
\end{prop}
 
\medskip

We shall set in the proof 
$$
\mathbb{D}_n := \big\{r^n_k = k2^{-n}T\,;k=0\ldots 2^n\big\},
$$ 
for any $n\geq 1$, and talk about an element in one of the sets $\mathbb{D}_n$ as a dyadic times. Let insist here on the convention that $L^a$ stands for the integrability class of random variables, whereas we shall always write $L^a(\RR^d)$ for integrable functions on $\RR^d$.

\medskip
 
\begin{DemThm}
We first show that $V$ has a modification that is almost surely $\frac{1}{p'}$-H\"older, with values in $\cC^{2+r-\frac{d}{a}}(\RR^d,\RR^d)$. This is done in an elementary way using Besov embedding theorem to write
\begin{align*}
\EE\Big[&\|V_{t,s}\|^{2a}_{\cC^{2+r-\frac{d}{a}}}\Big]^\frac{1}{2a} \leq \EE\Big[\big\|V_{ts}\big\|^{2a}_{B_{a,a}^{r+2}}\Big]^\frac{1}{2a}   \\
& \lesssim  \EE\Big[\|V_{ts}\|_{L^a(\RR^d)}^{2a}\Big]^\frac{1}{2a} + \EE\left[\left(\int_{B(0,1)} \big\| \Delta^{k_1}_\sigma V_{ts}\big\|_{L^a(\RR^d)}^a \, \frac{d\sigma} {|\sigma|^{(2+r)a + d}} \right)^2\right]^\frac{1}{2a}   \\
& \lesssim  \left(\int_{\RR^d}\EE\left[|V_{ts}(y)|^{2a}\right]^\frac{1}{2} \, dy\right)^\frac{1}{a} + \left(\int_{\RR^d}\int_{B(0,1)} \EE\left[\big|\Delta^{k_1}_\sigma V_{ts}(y)\big|^{2a}\right]^\frac{1}{2}\, \frac{d\sigma} {|\sigma|^{(2+r)a + d}}\,dy\right)^\frac{1}{a}   \\
& \lesssim |t-s|^\frac{1}{p},
\end{align*}
with a multiplicative constant in the inequality proportional to $\big\|C_0^V\big\|+ \big\|C_1^V\big\|$. The result for $V$ follows then from the usual Kolmogorov regularity theorem, here for a process with values in $\cC^{2+r-\frac{d}{a}}(\RR^d,\RR^d)$; write $\widetilde{V}$ for its modification with values in $\cC^\frac{1}{p'}\big(D_T;\cC^{2+r-\frac{d}{a}}(\RR^d,\RR^d)\big)$. 

\ssk

Note, for $s\leq u\leq t$, the elementary inequality
\begin{equation}
\label{EqEstimateBracket}
\begin{split}
\Big\|\big[V_{us}, V_{tu}\big]\Big\|_{\cC^{1+r-\frac{d}{a}}} &\lesssim \big\| DV_{tu}\big\|_{\cC^{1+r-\frac{d}{a}}} \|V_{us}\|_{\cC^{2+r-\frac{d}{a}}} +\|DV_{u,s}\|_{\cC^{1+r-\frac{d}{a}}} \|V_{t,u}\|_{\cC^{2+r-\frac{d}{a}}}  \\
&\lesssim |t-u|^\frac{1}{p'}\,|u-s|^\frac{1}{p'}\,\big\|\widetilde{V}\big\|^2,
\end{split}
\end{equation}
where $\|\widetilde{V}\big\|$ stands for the norm $\widetilde{V}$ as an element of $\cC^\frac{1}{p'}\big(D_T;\cC^{2+r-\frac{d}{a}}(\RR^d,\RR^d)\big)$, and is in $L^a$ -- Proposition \ref{prop:product} is used to justify the first inequality.

\medskip

Given two dyadic times $s<t$, with $s=k_s2^{-n_0}$ and $t=k_t2^{-n_0}$, the interval $[s,t)$ can uniquely be written as a finite disjoint union of intervals $[u,v)$ with ends in $\mathbb{D}_n$, for $n\geq n_0+1$, and where no three intervals have the same length. Write $s=s_0 < s_1< \cdots<s_N < s_{N+1} = t$, for the induced partition of $[s,t)$, and note that 
\begin{equation}
\label{EqEstimateIncrements}
\sum_{n=0}^N \big(s_{n+1}-s_n\big)^\frac{1}{p'}\big(s_n - s_0\big)^\frac{1}{p'} \lesssim c\,|t-s|^\frac{2}{p'},
\end{equation}
for an absolute positive constant $c$. Using repeatedly the decomposition 
$$
W_{s_ns_0} = W_{s_ns_{n-1}} + W_{s_{n-1}s_0} + \frac{1}{2}\,\Big[V_{s_{n-1}s_0}\, , \,V_{s_ns_{n-1}}\Big],
$$
together with estimate \eqref{EqEstimateBracket} and \eqref{EqEstimateIncrements}, we see that 
$$
\big\| W_{ts}\big\|_{\cC^{1+r'-\frac{d}{a}}} \lesssim\Big(\|\widetilde{V}\big\|+ M\Big)\,|t-s|^\frac{2}{p'},
$$
where 
$$
M := \sum_{n\geq 0}\frac{2^{-n\big(1+\frac{2a}{p}\big)}}{(n+1)^2}\,\sum_{k=0}^{2^n-1} \Big\|W_{r_{k+1}^n r_k^n}\Big\|_{\cC^{1+r-\frac{d}{a}}}
$$
is an integrale random variable, so is almost surely finite, as a consequence of Besov embedding and assumptions \eqref{eq:kolmo1} and \eqref{eq:kolmog} on the vector field $W$. An obvious extension procedure, such as classically done in the proof of Kolmogorov regularity theorem, finishes the proof of the statement.
\end{DemThm}

\medskip 

\begin{thm}[Kolmogorov-Lamperti-type tightness criterion for rough drivers]   {\it 
\label{corollary:lamperti}
Let $\big(V^\epsilon,W^\epsilon\big)$ be a family of vector fields satisfying the assumptions of theorem \ref{theorem:kolmo}, with 
$$
\Big\|C_0^{V^\epsilon}\Big\|+ \Big\|C_0^{W^\epsilon}\Big\| + \Big\|C_1^{V^\epsilon}\Big\| + \Big\|C_1^{W^\epsilon}\Big\|
$$
uniformly bounded above as $\epsilon$ ranges in $(0,1]$. Then, for every positive $p',r'$ with $p'< 2 + r'$, and
$$
r'<r-\frac{d}{a},\qquad \textrm{ and }\qquad  \frac{1}{3}<\frac{1}{p'}<\frac{1}{p}-\frac{1}{2a},
$$
the family $\widetilde{V}^\epsilon$ is tight in the space of rough drivers with regularity indices $p'$ and $(2+r')$.   }
\end{thm}
  
\medskip  
 
\begin{Dem}
The proof is elementary and consists in using first theorem \ref{theorem:kolmo} with $p''>p'$ and $r''>r'$ satisfying the conditions, and seeing that the quantities 
$$
\EE\left[\big\|\widetilde{\bV}\big\|_{(p',2+r')}^{2p}\right]
$$ 
are bounded uniformly in $0<\epsilon\leq 1$. So the probability that $\widetilde{V}^\epsilon$ is outside a fixed ball in the space of $(p'',2+r'')$-rough drivers can be made arbitrarily small by choosing a large enough radius for that ball. The claim follows from the fact that such a ball is compact in the space of $(p',2+r')$-rough drivers, by a standard Ascoli-Arzela-type argument.
\end{Dem}

\medskip

Note that one can find in \cite{BailleulRiedel} other regularization and compactness results for rough drivers -- Theorem 28. The present Lamperti-type compactness criterion happens to be particularly easy to use in our study of stochastic turbulence, in section \ref{SectionKestenPapanicolaou}.

\bigskip

\section[\hspace{0.6cm} {\sf An elementary lemma}]{An elementary lemma}
\label{SectionAppendixLemma}

We state and prove here the following elementary lemma that was used in section \ref{SubsubsectionFiniteFimCvgce} to bring back the study of the convergence problem for a \textit{continuous family} of processes to the convergence problem for a \textit{sequence} of processes. We adopt here the notations of that section.

\medskip 

\begin{lem}\label{LemAppendix}   {\it 
If the finite dimensional laws of the sequence of processes 
$$
\Big(\fV^{\epsilon_n}_\bullet,\int_0^\bullet \fV^{\epsilon_n}_{\dd u}\star \fV^{\epsilon_n}_u\Big)(\xx)
$$ 
converge to some limit process then the finite dimensional laws of the continuous family 
$$
\Big(\fV^{\eps}_\bullet,\int_0^\bullet \fV^{\eps}_{\dd u}\star \fV^{\eps}_u\Big)(\xx),
$$ 
indexed by $0<\epsilon\leq 1$, also converge to the same limit.   }
\end{lem} 

\medskip

\begin{Dem}
We use for the first component the same argument as in the proof of proposition \ref{lemma:discretization2}, and use the fact that there is almost surely an exponent $b<1$ such that the processes $\frak{V}^\epsilon_\bullet$ have finite $b$-H\"older norm, uniformly in $\epsilon$. So, taking $n=\big[\eps^{-2}\big] $, we have the almost surely estimate
\begin{align*}
\Big|\fV^{\eps}_t(\xx) -\fV^\frac{1}{n}_t(\xx)\Big| &\leq \left|  \eps \int_{nt}^{\eps^{-2}t}\fF_r(\xx)\right| + \left| (\eps-n^{-\frac{1}{2}}) \int_{0}^{nt}\fF_r(\xx)\right|   \\
&\leq \big|\fV^{\eps}_{t}(\xx) -\fV^{\eps}_{t n\eps^2}(\xx) \big| + \left| (\eps n^{\frac{1}{2}}-1) \fV_t^\frac{1}{n}(\xx)\right|  \\
&\lesssim t^\beta\big(|1-n\eps^2|^\beta + |\eps n^{\frac{1}{2}}-1| \big)\lesssim \big(\eps^{2\beta}+ \eps^2\big).
\end{align*} 
The proof for the second component is similar. Write
\begin{align*}
\left|\int_0^t \fV^{\eps}_{du} \right. & \left.\star \fV^{\eps}_u(\xx) - \int_0^t \fV^\frac{1}{n}_{du} \star \fV^\frac{1}{n}_u(\xx)\right|   \\
&\leq \left| (1-\eps^{-2}n^{-1}) \int_0^{n\eps^2 t} \fV^{\eps}_{du} \star \fV^{\eps}_u(\xx)\right|   \\
&+ \left|\eps^2 \int_{nt}^{\eps^{-2}t} \int_{nt}^{u_1} \fF_{u_1} \star \fF_{u_2}(\xx) \,du_2du_1\right|  + \eps^2 n \left|\Big(\int_{nt}^{\eps^2t} \fF_u\,du\Big)\star \Big(n^{-1}\int_0^{nt} \fF_u \dd u\Big)\right|.
\end{align*}
Again, by Theorem \ref{theorem:tighness} and by the definition of $n$, the first term of the right hand side is bounded by an almost surely finite constant multiple of $\eps^2$. Since, $\fF_u$ is almost surely bounded, the second term of the right hand side is of order $\eps^2$. Since $\eps^2 n \int_{nt}^{\eps^2 t} \fF_u(\xx)\,du$ is almost surely bounded by a constant independent of $\epsilon$, we eventually have the estimate 
$$
\Big|\int_0^t \fV^{\eps}_{du} \star \fV^{\eps}_u(\xx) - \int_0^t \fV^{\epsilon_n}_{du} \star \fV^{\epsilon_n}_u(\xx)\Big| \lesssim \eps^2 + \left|\frac1n\int_0^{nt} \fF_u(\xx) \,du\right|.
$$
The conclusion follows from the fact that $\fF_\bullet(\xx)$ is centered, stationary and mixing, from which the ergodic theorem implies that $\frac1n\int_0^{nt} \fF_u(\xx) \, du$ tends to $0$. 
\end{Dem}


\vfill

\end{document}